\documentclass[12pt]{amsart}
\usepackage{amssymb}
\begin{document}

\def\dbl{[\hskip -1pt[}
\def\dbr{]\hskip -1pt]}

\title{Approximation and convergence of formal CR-mappings}
\author[F. Meylan, N. Mir, D. Zaitsev]{Francine Meylan, Nordine Mir, and Dmitri Zaitsev}
\address{F. Meylan: Institut de Math\'ematiques, Universit\'e de Fribourg, 1700 Perolles, Fribourg, Switzerland}
\email{francine.meylan@unifr.ch}
\address{N. Mir: Universit\'e de Rouen, Laboratoire de Math\'ematiques Rapha\"el Salem, UMR 6085 CNRS, 76821 Mont-Saint-Aignan Cedex, France}
\email{Nordine.Mir@univ-rouen.fr}
\address{D. Zaitsev: Mathematisches Institut, Eberhard-Karls-Universit\"at T\"ubingen, Auf der Morgenstelle 10,
72076 T\"ubingen, Germany}
\email{dmitri.zaitsev@uni-tuebingen.de}
\thanks{\noindent 2000 {{\em Mathematics Subject Classification.} 32H02, 32V20, 32V40 \\
The first author was partially supported by Swiss NSF Grant 2100-063464.00/1}}

\def\Label#1{\label{#1}}
\def\1#1{\ov{#1}}
\def\2#1{\widetilde{#1}}
\def\6#1{\mathcal{#1}}
\def\4#1{\mathbb{#1}}
\def\3#1{\widehat{#1}}
\def\K{{\4K}}
\def\LL{{\4L}}

\def\C{{\4C}}
\def\R{{\4R}}
\def \MM{{\4M}}

\def\Re{{\sf Re}\,}
\def\Im{{\sf Im}\,}

\numberwithin{equation}{section}
\def\s{s}
\def\k{\kappa}
\def\ov{\overline}
\def\span{\text{\rm span}}
\def\ad{\text{\rm ad }}
\def\tr{\text{\rm tr}}
\def\xo {{x_0}}
\def\Rk{\text{\rm Rk\,}}
\def\sg{\sigma}
\def \emxy{E_{(M,M')}(X,Y)}
\def \semxy{\scrE_{(M,M')}(X,Y)}
\def \jkxy {J^k(X,Y)}
\def \gkxy {G^k(X,Y)}
\def \exy {E(X,Y)}
\def \sexy{\scrE(X,Y)}
\def \hn {holomorphically nondegenerate}
\def\hyp{hypersurface}
\def\prt#1{{\partial \over\partial #1}}
\def\det{{\text{\rm det}}}
\def\wob{{w\over B(z)}}
\def\co{\chi_1}
\def\po{p_0}
\def\fb {\bar f}
\def\gb {\bar g}
\def\Fb {\ov F}
\def\Gb {\ov G}
\def\Hb {\ov H}
\def\zb {\bar z}
\def\wb {\bar w}
\def \qb {\bar Q}
\def \t {\tau}
\def\z{\chi}
\def\w{\tau}
\def\Z{\zeta}
\def\phi{\varphi}
\def\eps{\varepsilon}

\def \T {\theta}
\def \Th {\Theta}
\def \L {\Lambda}
\def\b {\beta}
\def\a {\alpha}
\def\o {\omega}
\def\l {\lambda}

\def \im{\text{\rm Im }}
\def \re{\text{\rm Re }}
\def \Char{\text{\rm Char }}
\def \supp{\text{\rm supp }}
\def \codim{\text{\rm codim }}
\def \Ht{\text{\rm ht }}
\def \Dt{\text{\rm dt }}
\def \hO{\widehat{\mathcal O}}
\def \cl{\text{\rm cl }}
\def \bR{\mathbb R}
\def \bS{\mathbb S}
\def \bK{\mathbb K}
\def \bD{\mathbb D}
\def \bC{\mathbb C}
\def \C{\mathbb C}
\def \N{\mathbb N}
\def \bL{\mathbb L}
\def \bZ{\mathbb Z}
\def \bN{\mathbb N}
\def \scrF{\mathcal F}
\def \scrK{\mathcal K}
\def \mc #1 {\mathcal {#1}}
\def \scrM{\mathcal M}
\def \cR{\mathcal R}
\def \scrJ{\mathcal J}
\def \scrA{\mathcal A}
\def \scrO{\mathcal O}
\def \scrV{\mathcal V}
\def \scrL{\mathcal L}
\def \scrE{\mathcal E}
\def \hol{\text{\rm hol}}
\def \aut{\text{\rm aut}}
\def \Aut{\text{\rm Aut}}
\def \J{\text{\rm Jac}}
\def\jet#1#2{J^{#1}_{#2}}
\def\gp#1{G^{#1}}
\def\gpo{\gp {2k_0}_0}
\def\emmp {\scrF(M,p;M',p')}
\def\rk{\text{\rm rk\,}}
\def\Orb{\text{\rm Orb\,}}
\def\Exp{\text{\rm Exp\,}}
\def\Span{\text{\rm span\,}}
\def\d{\partial}
\def\D{\3J}
\def\pr{{\rm pr}}

\def \CZZ {\C \dbl Z,\zeta \dbr}
\def \D{\text{\rm Der}\,}
\def \Rk{\text{\rm Rk}\,}
\def \CR{\text{\rm CR}}
\def \ima{\text{\rm im}\,}
\def \I {\mathcal I}
\def \M {\mathcal M}

\newtheorem{Thm}{Theorem}[section]
\newtheorem{Cor}[Thm]{Corollary}
\newtheorem{Pro}[Thm]{Proposition}
\newtheorem{Lem}[Thm]{Lemma}

\theoremstyle{definition}\newtheorem{Def}[Thm]{Definition}

\theoremstyle{remark}
\newtheorem{Rem}[Thm]{Remark}
\newtheorem{Exa}[Thm]{Example}
\newtheorem{Exs}[Thm]{Examples}

\def\bl{\begin{Lem}}
\def\el{\end{Lem}}
\def\bp{\begin{Pro}}
\def\ep{\end{Pro}}
\def\bt{\begin{Thm}}
\def\et{\end{Thm}}
\def\bc{\begin{Cor}}
\def\ec{\end{Cor}}
\def\bd{\begin{Def}}
\def\ed{\end{Def}}
\def\br{\begin{Rem}}
\def\er{\end{Rem}}
\def\be{\begin{Exa}}
\def\ee{\end{Exa}}
\def\bpf{\begin{proof}}
\def\epf{\end{proof}}
\def\ben{\begin{enumerate}}
\def\een{\end{enumerate}}

\keywords{Formal holomorphic map, convergence, real-analytic and real-algebraic
CR-submanifolds, Artin's approximation Theorem}

\maketitle

\begin{abstract} Let $M\subset \C^N$ be a minimal real-analytic
CR-submanifold and $M'\subset \C^{N'}$ a real-algebraic subset
through points $p\in M$ and $p'\in M'$. We show that that any
formal (holomorphic) mapping $f\colon (\C^N,p)\to (\C^{N'},p')$,
sending $M$ into $M'$, can be approximated up to any given order
at $p$ by a convergent map sending $M$ into $M'$. If $M$ is
furthermore generic, we also show that any such map $f$, that is
not convergent, must send (in an appropriate sense) $M$ into the
set $\6E'\subset M'$ of points of {\sc D'Angelo} infinite type.
Therefore, if $M'$ does not contain any nontrivial
complex-analytic subvariety through  $p'$, any formal map $f$ as
above is necessarily convergent.
\end{abstract}


\section{Introduction and results}\Label{int}
An important step in understanding the existence of analytic
objects with certain properties consists of understanding the same
problem at the level of formal power series. The latter problem
can be reduced to a sequence of algebraic equations for the
coefficients of the unknown power series and is often simpler than
the original problem, where the power series are required to be
convergent. It is therefore of interest to know whether such power
series are automatically convergent or can possibly be replaced by
other convergent power series satisfying the same properties. A
celebrated result of this kind is {\sc Artin}'s Approximation
Theorem \cite{A68} which states that a formal solution of a system
of analytic equations can be replaced by a convergent solution of
the same system that approximates the original solution at any
prescribed order.

In this paper we study convergence and approximation properties
(in the spirit of \cite{A68}) of formal (holomorphic) mappings sending
real-analytic submanifolds $M\subset\C^N$ and $M'\subset\C^{N'}$ into each other.
In this situation the above theorem of {\sc Artin} cannot be applied directly.
Moreover, without additional assumptions on the submanifolds,
the analogous approximation statement is not even true.
Indeed,
in view of an example of {\sc Moser-Webster} \cite{MW},
there exist real-algebraic surfaces  $M,M'\subset\C^2$
that are formally but not biholomorphically equivalent.
However, our first main result shows that this phenomenon cannot happen
if $M$ is a minimal CR-submanifold (not necessarily algebraic) in $\C^N$
(see \S\ref{sarra} for the notation and definitions):

\begin{Thm}\Label{approx}
Let $M\subset \C^N$ be a real-analytic minimal CR-submanifold
and $M'\subset \C^{N'}$ a real-algebraic subset
with $p\in M$ and $p'\in M'$.
Then for any formal $($holomorphic$)$ mapping $f\colon (\C^N,p)\to (\C^{N'},p')$ sending
$M$ into $M'$ and any positive integer $k$, there exists a germ of
a holomorphic map $f^{k}\colon (\C^N,p)\to (\C^{N'},p')$ sending $M$ into $M'$,
whose Taylor series at $p$ agrees with $f$ up to order $k$.
\end{Thm}

Approximation results in the spirit of Theorem~\ref{approx} have
been recently obtained in \cite{BRZ,BMR} in the important case
when $N=N'$ and $f$ is invertible. Note that under the assumptions
 of Theorem \ref{approx}, there may exist
 nonconvergent maps $f$ sending $M$ into $M'$.
For instance, it is easy to construct such maps in case $M$ is not
generic in $\C^N$. Also, if $M'$ contains an irreducible
complex-analytic subvariety $E'$ of positive dimension through
$p'$, such maps $f$ with $f(M)\subset E'$ (in the formal sense)
always exist. Our next result shows that these are essentially the
only exceptions. Denote by $\6 E'$ the set  of all points of $M'$
through which there exist irreducible complex-analytic
subvarieties of $M'$ of positive dimension. This set is always
closed (\cite{Le86,D}) but not real-analytic in general (see
\cite{MMZ01} for an example). In the following, we say that a
formal (holomorphic) map $f\colon (\C^N,p)\to (\C^{N'},p')$ sends
$M$ into $\6E'$ if $\phi(f(x(t)))\equiv 0$ holds for all germs of
real-analytic maps $x\colon (\R^{\dim M}_t,0)\to (M,p)$ and
$\phi\colon(M',p')\to(\R,0)$ such that $\phi$ vanishes on $\6 E'$.
We prove:

\begin{Thm}\Label{ncc0}
Let $M\subset \C^N$ be a minimal real-analytic generic submanifold
and $M'\subset \C^{N'}$ a real-algebraic subset
with $p\in M$ and $p'\in M'$.
Then any formal $($holomorphic$)$ mapping $f\colon (\C^N,p)\to (\C^{N'},p')$ sending
$M$ into $M'$ is either convergent or sends $M$ into $\6E'$.
\end{Thm}

As an immediate consequence we obtain the following
characterization:

\begin{Cor}\Label{stock}
Let $M\subset \C^N$ be a minimal real-analytic generic submanifold
and $M'\subset \C^{N'}$ a real-algebraic subset with $p\in M$ and $p'\in M'$.
Then all formal maps $f\colon (\C^N,p)\to (\C^{N'},p')$ sending $M$ into $M'$
are convergent if and only if $M'$ does not contain
any irreducible complex-analytic subvariety of positive dimension through $p'$.
\end{Cor}

In contrast to most previously known related results,
Theorems~\ref{approx}--\ref{ncc0} and Corollary~\ref{stock} do not
contain any assumption on the map $f$. Indeed,
Theorems~\ref{approx}--\ref{ncc0} seem to be the first results of
this kind and an analog of Corollary~\ref{stock} appears only in
the work of {\sc Baouendi-Ebenfelt-Rothschild} \cite{BER00} for
the case $M,M'\subset\C^N$ are real-analytic hypersurfaces
containing no nontrivial complex subvarieties. In fact they prove
a more general result for $M$ and $M'$ of higher codimension
assuming the map $f$ to be finite and show (see the proof of
\cite[Proposition~7.1]{BER00}) that the finiteness of $f$
automatically holds (unless $f$ is constant) in the mentioned case
of hypersurfaces. However, in the setting of
Corollary~\ref{stock}, the finiteness of a (nonconstant) map $f$
may fail to hold even when $M,M'\subset\C^N$ are hypersurfaces,
e.g. for $M:=S^3\times\C, M':=S^5\subset\C^3$, where
$S^{2n-1}\subset\C^n$ is the unit sphere. Thus, even in this case,
Theorems~\ref{approx}--\ref{ncc0} and Corollary~\ref{stock} are
new and do not follow from the same approach. It is worth
mentioning that Corollary~\ref{stock} is also new in the case of
unit spheres $M=S^{2N-1}$ and $M'=S^{2N'-1}$ with $N'>N$.

Previous work in the direction of Theorem~\ref{ncc0} is due to
{\sc Chern-Moser} \cite{CM} for real-analytic Levi-nondegenerate
hypersurfaces. More recently, this result was extended in
\cite{BER97, BER99, BER00, BRZ, M00, M01, BMR, la} under weaker
conditions on the submanifolds and mappings.

One of the main novelties of this paper compared to previous
related work lies in the study of convergence properties of ratios
of formal power series rather than of the series themselves. It is
natural to call such a ratio convergent if it is equivalent to a
ratio of convergent power series. However, for our purposes, we
need a refined version of convergence along a given submanifold
that we define in \S\ref{warm-up} (see
Definition~\ref{convergent}). With this refined notion, we are
able to conclude the convergence of a given ratio along a
submanifold provided its convergence is known to hold along a
smaller submanifold and  under suitable conditions on the ratio
(see Lemmata~\ref{trivial}--\ref{nontrivial}).

Another novelty of our techniques consists of applying the
mentioned convergence results of \S\ref{warm-up} and their
consequences given in \S\ref{pullback} to ratios defined on
iterated complexifications of real-analytic submanifolds (in the
sense of \cite{Z97,Z99}) rather than on single Segre sets (in the
sense of \cite{BERacta}) associated to given fixed points. The
choice of iterated complexifications is needed to guarantee the
nonvanishing of the relevant ratios that may not hold when
restricted to the Segre sets. These tools are then used to obtain
the convergence of a certain type of ratios of formal power series
that appear naturally in the proofs of
Theorems~\ref{approx}--\ref{ncc0}. This is done in
Theorem~\ref{hurry} that is, in turn, derived from
Theorem~\ref{mer-ext} which is established in the more general
context of a pair of submersions of a complex manifold.

After the necessary preparations in \S\S\ref{Zar}--\ref{NEW1}, we
state and prove Theorem~\ref{straight} which is the main technical
result of the paper and which implies, in particular, that the
(formal) graph of $f$ is contained in a real-analytic subset
$Z_f\subset M\times M'$ satisfying a straightening property. If
$f$ is not convergent, the straightening property implies the
existence of nontrivial complex-analytic subvarieties in $M'$ and
hence proves Theorem~\ref{ncc0}. To obtain Theorem~\ref{approx},
we use the additional property of the set $Z_f$ (also given by
Theorem~\ref{straight}), stating that $Z_f$ also contains graphs
of holomorphic maps approximating $f$ up to any order (at $0$).
The fact that $Z_f\subset M\times M'$ then yields
Theorem~\ref{approx}.

\section{Notation and definitions}

\subsection{Formal mappings and CR-manifolds}\Label{sarra}
A formal (holomorphic) mapping $f:(\C_Z^N,p)\to (\C_{Z'}^{N'},p')$
is the data of $N'$ formal power series $(f_1,\ldots,f_{N'})$ in
$Z-p$, with $f(p)=p'$. Let $M\subset \C^N$ and $M'\subset \C^{N'}$
be real-analytic submanifolds with $p\in M$ and $p'\in M'$, and
$\rho (Z,\1Z)$, $\rho'(Z',\1{Z'})$ be real-analytic vector-valued
defining functions for $M$ near $p$ and $M'$ near $p'$
respectively. Recall that a formal mapping $f$ as above sends $M$
into $M'$ if there exists a matrix $a(Z,\1Z)$, with entries in $\C
\dbl Z-p,\1Z-\1p\dbr$, such that the formal identity
\begin{equation}\Label{mafia}
\rho'(f(Z),\1{f(Z)})=a(Z,\1Z)\cdot \rho (Z,\1Z)
\end{equation}
holds. Observe that (\ref{mafia}) is independent of the choice of
local real-analytic defining functions for $M$ and $M'$. For $M'$
merely a real-analytic subset in $\C^{N'}$, we also say that $f$
sends $M$ into $M'$, and write $f(M)\subset M'$, if (\ref{mafia})
holds for any real-analytic function $\rho'$ (with some $a$
depending on $\rho'$), defined in a neighborhood of $p'$ in
$\C^{N'}$, vanishing on $M'$. The notation $f(M)\subset M'$ is
motivated by the fact that in case $f$ is convergent, the above
condition holds if and only if $f$ is the Taylor series of a
holomorphic map sending $(M,p)$ into $(M',p')$ in the sense of
germs.

For a real-analytic CR-submanifold $M\subset \C^N$ (see e.g.\
\cite{BERbook} for basic concepts related to CR-geometry), we
write $T_p^cM$ for the complex tangent space of $M$ at $p\in M$,
i.e.\ $T_p^cM:=T_pM\cap i T_pM$. Recall that $M$ is called {\em
generic} if for any point $p\in M$, one has $T_pM+iT_pM=T_p\C^N$.
Recall also that $M$ is called {\em minimal} (in the sense of {\sc
Tumanov} \cite{T}) at a point $p\in M$ if there is no real
submanifold $S\subset M$ through $p$ with ${\rm dim}\, S<{\rm
dim}\, M$ and $T_q^cM\subset T_qS$ for all $q\in S$. It is
well-known that, if $M$ is real-analytic, the minimality of $M$ at
$p$ is equivalent to the finite type condition of {\sc Kohn}
\cite{Kohn} and {\sc Bloom-Graham} \cite{BG}.

\subsection{Rings of formal power series}\Label{formalstuff}
For a positive integer $n$, we  write $\C\dbl t\dbr$ for the ring
of formal power series (with complex coefficients) in the
indeterminates $t=(t_1,\ldots,t_n)$ and $\C \{t\}$ for the ring of
convergent ones. If $t^{0}\in \C^n$, $\C\dbl t-t^0 \dbr$ and $\C
\{t-t^0\}$ will denote the corresponding rings of series centered
at $t^0$. For any formal power series $F(t)$, we denote by $\1{F}
(t)$ the formal power series obtained from $F(t)$ by taking
complex conjugates of its coefficients.

An ideal $I\subset \C \dbl t\dbr$ is called a {\em manifold ideal}
if it has a set of generators with linearly independent
differentials (at $0$). If $I\subset \C \dbl t\dbr$ is a manifold
ideal, then any set of generators with linearly independent
differentials has the same number of elements that we call the
{\em codimension} of $I$. In general, we say that a manifold ideal
$I$ defines a {\em formal submanifold} $\6{S}\subset \C^l$ and write
$I=I(\6 S)$. Note that if $I\subset \C \{t\}$, then $I$ defines a
(germ of a) complex submanifold $\6 S\subset \C^n$ through the origin in the
usual sense. Given a formal submanifold $\6S\subset \C^n$ of
codimension $d$, a (local) parametrization of $\6S$ is a formal map
$j\colon (\C^{n-d},0)\to (\C^n,0)$ of rank $n-d$ (at $0$)
such that $V \circ j=0$ for all $V\in I(\6{S})$. If
$\6{S},\6{S}'\subset \C^n$ are two formal submanifolds, we write
$\6{S}\subset \6{S}'$ to mean that $I(\6{S}')\subset I(\6{S})$.
For a formal map $h:(\C_t^n,0)\to (\C_{T}^r,0)$, we define its
graph $\Gamma_h\subset \C^n\times \C^r$ as the
 formal submanifold given by $I(\Gamma_h)$, where $I(\Gamma_h)\subset \C \dbl t,T\dbr$ is the ideal
 generated by $T_1-h_1(t),\ldots,T_r-h_r(t)$.

For a formal power series $F(t)\in \C \dbl t\dbr$ and a formal
submanifold $\6S\subset \C^l$, we write $F|_{\6 S}\equiv 0$ (or
sometimes also $F(t)\equiv 0$ for $t\in \6S$) to mean that
$F(t)\in I(\6{S})$. If $k$ is a nonnegative integer, we also write
$F(t)=O(k)$ for $t\in \6{S}$ to mean that for one (and hence for
any) parametrization $j=j(t)$ of $\6{S}$, $(F\circ j)(t)$ vanishes
up to order $k$ at the origin. We also say that another power series $G(t)$
agrees with $F(t)$ up to order $k$ (at the origin) if $F(t)-G(t)=O(k)$.

A convenient criterion for the convergence of a formal power
series is given by the following well-known result (see e.g.\
\cite{BER00,M00} for a proof).

\bp\Label{recall}
Any formal power series which satisfies a
nontrivial polynomial identity with convergent coefficients is
convergent.
\ep

It will be also convenient to consider formal power series
defined on an abstract complex manifold (of finite dimension) $\6X$
centered at a point $x_0\in \6X$ without referring to specific coordinates.
In each coordinate chart such a power series
is given by a usual formal power series
that transforms in the obvious way under biholomorphic coordinate changes.
Given such a series $H$, we write $H(x_0)$ for the value at $x_0$ that is always defined.
It is easy to see that the set of all
formal power series on a complex manifold
centered at $x_0$ forms a (local) commutative ring that is an integral domain.
The notion of convergent power series extends to
power series on abstract complex manifolds in the obvious way.

In a similar way, one may consider formal holomorphic vector
fields on  abstract complex manifolds and apply them to formal
power series. If $F$ and $G$ are such formal power series on $\6
X$ centered at $x_0$, we write $\6 L (F/G)\equiv 0$ if and only if
$F \6L G-G\6L F\equiv 0$ (as formal power series on $\6 X$).

Completely analogously one may define formal power series mappings
between complex manifolds and their compositions.

\section{Meromorphic extension of ratios of formal power series}\Label{fou}

The ultimate goal of this section  is to establish a meromorphic
extension property for ratios of formal power series (see Theorem
\ref{mer-ext}).

\subsection{Convergence of ratios of formal power series}\Label{warm-up}

Throughout \S \ref{fou}, for any formal power series $F=F(t)\in\C
\dbl t\dbr$ in $t=(t_1,\ldots,t_n)$ and any nonnegative integer
$k$, we denote by $j^k F$ or by $j_t^kF$ the formal power series
mapping corresponding to the collection of all partial derivatives
of $F$ up to order $k$. We shall use the first notation when there
is no risk of confusion and the second one when other
indeterminates appear. For $F(t),G(t)\in \C\dbl t\dbr$, we write
$(F:G)$ for a pair of two formal power series thinking of it as a
ratio, where we allow both series to be zero.

\begin{Def}\Label{k-similar}
Let $(F_1:G_1)$, $(F_2:G_2)$ be ratios of formal power series in
$t=(t_1,\ldots,t_n)$, and $S\subset\C^n$ be a (germ of a) complex
submanifold through $0\in \C^n$. We say that the ratios
$(F_1:G_1)$ and $(F_2:G_2)$ are {\em $k$-similar along $S$} if
$(j^k(F_1 G_2 - F_2 G_1))|_S \equiv 0$.
\end{Def}
The defined relation of similarity for formal power series is
obviously symmetric but not transitive, e.g.\ any ratio is
$k$-similar to $(0:0)$ along any complex submanifold $S$ and for
any nonnegative integer $k$. However, we have the following weaker
property:

\bl\Label{transitivity} Let $(F_1:G_1)$, $(F_2:G_2)$  and
$(F_3:G_3)$ be ratios of formal power series in
$t=(t_1,\ldots,t_n)$, $S\subset \C^n$ a complex submanifold
through the origin and $k$ a nonnegative integer. Suppose that
both ratios $(F_1:G_1)$ and $(F_3:G_3)$ are $k$-similar to
$(F_2:G_2)$ along $S$. Then, if there exists $l\le k$ such that
$(j^l(F_2,G_2))|_S\not\equiv 0$, then $(F_1:G_1)$ and $(F_3:G_3)$
are $(k-l)$-similar along $S$. \el

\bpf Without loss of generality, we may assume that $(j^l
F_2)|_S\not\equiv 0$. By the assumptions, we have $(j^k(F_1 G_2 - F_2
G_1))|_S \equiv 0$ and $(j^k(F_3 G_2 - F_2 G_3))|_S \equiv 0$.
Multiplying the first identity by $F_3$, the second by $F_1$ and
subtracting from each other, we obtain $(j^k(F_2(F_1 G_3 - F_3
G_1)))|_S \equiv 0$. Since $(j^l F_2)|_S\not\equiv 0$, the last
identity is only possible if $(j^{k-l}(F_1 G_3 - F_3 G_1))|_S
\equiv 0$ as required. \epf

We shall actually use the following refined version of
Lemma~\ref{transitivity} whose proof is completely analogous. In
what follows, for some splitting of indeterminates
$t=(t^1,t^2,t^3)\in \C^{n_1}\times \C^{n_2}\times
\C^{n_3}$ and for any formal power series $F(t)\in \C \dbl t\dbr$,
we write $j^k_{t^i}F$ for the collection of all partial
derivatives up to order $k$ of $F$ with respect to $t^i$,
$i=1,2,3$.

\bl\Label{fine-transitivity} Let $(F_1:G_1)$,  $(F_2:G_2)$ and
$(F_3:G_3)$ be ratios of formal power series in
$t=(t^1,t^2,t^3)\in\C^{n_1}\times \C^{n_2}\times
\C^{n_3}$ and set $S:=\C^{n_1}\times\{(0,0)\}\subset
\C^{n_1}\times \C^{n_2}\times \C^{n_3}$. Suppose that there  exist
integers $l\ge 0$, and $k_2,k_3\ge l$ such that
$(j^l(F_2,G_2))|_S\not\equiv 0$, $(j_{t^2}^{k_2}
j_{t^3}^{k_3} (F_1 G_2 - F_2 G_1))|_S \equiv 0$ and
$(j_{t^2}^{k_2} j_{t^3}^{k_3} (F_3 G_2 - F_2 G_3))|_S
\equiv 0$. Then $(j_{t^2}^{k_2-l} j_{t^3}^{k_3-l} (F_1 G_3
- F_3 G_1))|_S \equiv 0$. \el

Clearly, given a complex submanifold $S\subset \C^n$ through the
origin, any fixed ratios $(F_1:G_1)$ and $(F_2:G_2)$ are
$k$-similar along $S$ for any $k$ if and only
if $F_1G_2-F_2G_1\equiv 0$, i.e.\ if they are equivalent as ratios.
We now define a notion of convergence along $S$ for any ratio of
formal power series.

\bd\Label{convergent} Let $S\subset \C^n$ be a complex submanifold
through the origin and $F(t),G(t)\in \C \dbl t\dbr$,
$t=(t_1,\ldots,t_n)$. The ratio $(F:G)$ is said to be {\em
convergent along $S$} if there exist a nonnegative integer $l$
and, for any  nonnegative integer $k$,  convergent power series
$F_k(t),G_k(t)\in \C \{t\}$, such that the ratio $(F_k:G_k)$ is
$k$-similar to $(F:G)$ along $S$ and $(j^l(F_k,G_k))|_S\not\equiv 0$.
\ed

The uniformity of the choice of the integer $l$ is a crucial
requirement in Definition \ref{convergent} (see e.g.\ the proof of
Lemma \ref{nontrivial} below). This notion of convergence for
ratios of formal power series has the following elementary
properties.

\bl\Label{hue} For $F(t),G(t)\in \C \dbl t\dbr$, $t=(t_1,\ldots,t_n)$, the following hold. \ben
\item[(i)] $(F:G)$ is always convergent along $S=\{0\}$.
\item[(ii)] If $F$ and $G$ are convergent, then $(F:G)$ is convergent along any submanifold $S\subset\C^n$
$($through $0)$.
\item[(iii)] If $(F:G)$ is equivalent to a nontrivial ratio that is convergent along a submanifold $S$,
then $(F:G)$ is also convergent along $S$.
\item[(iv)] If $(F:G)$ is convergent along $S=\C^n$, then
it is equivalent to a nontrivial ratio of convergent power series.
\een \el

\begin{proof}
Properties (i), (ii) and (iv) are easy to derive from Definition
\ref{convergent}. Part (iii) is a consequence of Lemma \ref{transitivity}.
\end{proof}

An elementary useful property of ratios of formal power series is
given by the following lemma.

\bl\Label{talk} Let $(F:G)$ be a ratio of formal power series in
$t=(t^1,t^2)\in \C^{n_1}\times \C^{n_2}$ with $G\not
\equiv 0$, and such that $({\d}/{\d t^2})(F/G) \equiv 0$. Then there exists
$\2F(t^1),\2G(t^1)\in \C \dbl t^1\dbr$  with $\2G \not \equiv 0$
such that
$(F:G)$ is equivalent to $(\2F:\2G)$. \el

\begin{proof} From the assumption it is easy to obtain, by differentiation, the identity
\begin{equation}\Label{francine}
(\partial_{t^2}^{\nu}F)(\partial^{\alpha}_{t^2}G)-
(\partial^{\nu}_{t^2}G)(\partial^{\alpha}_{t^2}F)
 \equiv 0,
\end{equation}
for all multiindices $\alpha,\nu \in \N^{n_2}$. Since $G\not
\equiv 0$, there exists $\alpha \in \N^{n_2}$ such that
$(\partial^{\alpha}_{t^2}G)|_{t^2=0}\not \equiv 0$.
Define $\2F(t^1):=\partial^{\alpha}_{t^2}F(t^1,0)\in \C \dbl t^1\dbr $ and
$\2G(t^1):=\partial^{\alpha}_{t^2}G(t^1,0)\in \C \dbl t^1\dbr$. Then by putting $t^2=0$ in
(\ref{francine}), we obtain  that
$(\2F \partial^{\nu}_{t^2}G -
 \2G \partial_{t^2}^{\nu}F)|_{t^2=0} \equiv 0$ for any
multiindex $\nu \in \N^{n_2}$. From this, it follows that the
ratios $(F:G)$ and $(\2F:\2G)$ are equivalent, which completes the
proof of the lemma since by construction $\2G\not \equiv 0$.
\end{proof}

 The following lemma will be used in \S
\ref{pullback} to pass from smaller sets of convergence to larger
ones.

\bl\Label{trivial} Let $F(t),G(t)\in \C\dbl t \dbr$ be formal
power series in $t=(t^1,t^2,t^3)\in\C^{n_1}\times
\C^{n_2}\times \C^{n_3}$ that depend only on $(t^1,t^3)$.
Then, if the ratio $(F:G)$ is convergent along
$\C^{n_1}\times\{(0,0)\}$, it is also convergent along
$\C^{n_1}\times\C^{n_2}\times\{0\}$. \el

\begin{proof} We set $S:=\C^{n_1}\times\{(0,0)\}\subset \C^{n_1}\times
\C^{n_2}\times \C^{n_3}$ and $\2 S:=\C^{n_1}\times \C^{n_2}\times
\{0\} \subset \C^{n_1}\times \C^{n_2}\times \C^{n_3}$. By the
assumptions and Definition \ref{convergent}, there exists a
nonnegative integer $l$ and, for any nonnegative integer $k$,
convergent power series $F_k(t),G_k(t)\in \C \{t\}$ such that
\begin{equation}\Label{due}
\big(j_t^k\big(F(t^1,t^3)G_k(t)-G(t^1,t^3)F_k(t)\big)\big)|_S\equiv
0,
\end{equation}
and $(j^l(F_k,G_k))|_{S}\not \equiv 0$. We fix $k\geq l$. Choose
$\beta_0\in N^{n_2}$ with $|\beta_0|\leq l$ such that
$$\big(j^l_{(t^1,t^3)}\big((\partial^{\beta_0}_{t^2}
F_k)(t^1,0,t^3),(\partial^{\beta_0}_{t^2}
G_k)(t^1,0,t^3)\big)\big)|_{t^3=0}\not \equiv 0.$$ Define $\2
F_k(t):=\partial^{\beta_0}_{t^2}F_k(t^1,0,t^3)$ and
$\2 G_k(t):=\partial^{\beta_0}_{t^2}G_k(t^1,0,t^3)$.
By the construction, we have  $(j^l(\2F_k,\2G_k))|_{\2 S} \not \equiv
0$ and it is also easy to see from (\ref{due}) that $(F:G)$ is
$(k-l)$-similar to $(\2 F_k:\2 G_k)$ along $\2 S$. This finishes
the proof of Lemma \ref{trivial}.
\end{proof}

The next less obvious lemma will be also used for the same purpose.

\bl\Label{nontrivial} Consider formal power series $F(t),G(t)\in
\C\dbl t \dbr$ in $t=(t^1,t^2,t^3)\in\C^{n_1}\times
\C^{n_2}\times \C^{n_3}$ of the form
\begin{equation}\Label{nordine}
F(t)=\phi(Y(t^1,t^3),t^2), \quad
G(t)=\psi(Y(t^1,t^3),t^2),
\end{equation}
where $Y(t^1,t^3)\in (\C \dbl t^1,t^3 \dbr)^r$ for
some integer $r\ge 1$ and $\phi$ and $\psi$ are convergent power
series in $\C^r\times\C^{n_2}$ centered at $(Y(0),0)$. Then the
conclusion of Lemma~{\rm \ref{trivial}} also holds. \el

\bpf The statement obviously holds if $F$ and $G$ are both zero,
hence we may assume that $(F,G)\not\equiv 0$. Then there exists a
nonnegative integer $d$ such that $(j^d (F,G))|_S\not\equiv 0$,
where $S:=\C^{n_1}\times\{(0,0)\}\subset \C^{n_1}\times
\C^{n_2}\times \C^{n_3}$. Since $(F:G)$ is assumed to be convergent along $S$,
there exist a nonnegative integer $l$ and ratios $(F_s:G_s)$,
$s=0,1,\ldots$, of convergent power series such that $(F:G)$ is
$s$-similar to $(F_s:G_s)$ and
\begin{equation}\Label{forget}
(j^l(F_s,G_s))|_S\not \equiv 0
\end{equation}
for all $s$. Then, for any $k\ge l$ and $s\ge k+l$, we have
\begin{equation}\Label{jet-id}
(j^{s-k}_{t^2} j^k_{t^3} (F G_s - F_s G))|_S \equiv 0.
\end{equation}
In view of \eqref{nordine} we may rewrite \eqref{jet-id} in the
form
\begin{equation}\Label{conv-id}
R_{s,k}((j^k_{t^3} Y)(t^1,0),t^1) \equiv 0,
\end{equation}
where $R_{s,k}$ is a convergent power series in the corresponding
variables. We view \eqref{conv-id} as a system of analytic
equations $R_{s,k}(y,t^1)=0$ for $k\ge l$ fixed and $s\ge k+l$
arbitrary and $y(t^1):=(j^k_{t^3} Y)(t^1,0)$ as a
formal solution of the system. By applying {\sc Artin}'s
approximation theorem \cite{A68}, for any positive integer
$\kappa$, there exists a convergent solution $y^{\kappa}(t^1)$
agreeing up to order $\kappa$ (at $0\in \C^{n_1} $) with
$y(t^1)$ (and depending also on $k$) and satisfying
$R_{s,k}(y^{\kappa}(t^1),t^1)\equiv 0$ for all $s$ as above. It is
easy to see that there exists a convergent power series
$Y^{\kappa}(t^1,t^3)$ (e.g.\ a polynomial in $t^3$)
satisfying $(j^k_{t^3} Y^{\kappa})|_{t^3=0} \equiv
y^{\kappa}(t^1)$. Hence the power series
$\2F^{\kappa}_k(t):=\phi(Y^{\kappa}(t^1,t^3),t^2)$ and
$\2G^{\kappa}_k(t):=\psi(Y^{\kappa}(t^1,t^3),t^2)$ are
convergent and agree with $F(t)$ and $G(t)$ respectively up to
order $\kappa$. Therefore by choosing $\kappa$ sufficiently large
(depending on $k$), we may assume that
$(j^d(\2F^{\kappa}_k,\2G^{\kappa}_k))|_S \not\equiv 0$. In what
follows, we fix such a choice of $\kappa$. By our construction,
\eqref{conv-id} is satisfied with $Y$ replaced by $Y^{\kappa}$ and
thus  \eqref{jet-id} is satisfied with $(F,G)$ replaced by
$(\2F^{\kappa}_k,\2G^{\kappa}_k)$ i.e.\
\begin{equation}\Label{repere}
(j^{s-k}_{t^2} j^k_{t^3} (\2 F^{\kappa}_k G_s - F_s \2
G^{\kappa}_{k}))|_S \equiv 0.
\end{equation}
 In view of Lemma~\ref{fine-transitivity}, (\ref{jet-id}),
 (\ref{repere}) and (\ref{forget}) imply
\begin{equation}\Label{k-id}
(j^{s-k-l}_{t^2}
 j^{k-l}_{t^3} (F \2G^{\kappa}_k -
\2F^{\kappa}_k G))|_S \equiv 0.
\end{equation}
Since $s$ can be taken arbitrarily large, \eqref{k-id} implies
that $(F:G)$ and $(\2F^{\kappa}_k : \2G^{\kappa}_k)$ are
$(k-l)$-similar along $\2S:=\C^{n_1}\times\C^{n_2}\times\{0\}$.
Since $(j^d(\2F^{\kappa}_k,\2G^{\kappa}_k))|_S \not\equiv 0$
implies $(j^d(\2F^{\kappa}_k,\2G^{\kappa}_k))|_{\2S} \not\equiv
0$, the ratio $(F:G)$ is convergent along $\2S$ (in the sense of
Definition~\ref{convergent}) and the proof is complete. \epf

For the proof of Theorem \ref{mer-ext}, we shall also need the
following lemma.

\bl\Label{annamaria}
Let $\eta$ be a holomorphic map from a
neighborhood of $0$ in $\C^{r}$ into a neighborhood of $0$ in
$\C^n$, with $\eta (0)=0$, and $\alpha (t), \beta (t)\in \C \dbl
t\dbr$, $t=(t_1,\ldots,t_n)$. Suppose that there exists a $($germ
of a$)$ complex submanifold $S\subset \C^r$ through $0$
such that $\eta|_{S}:S \to \C^n$ has maximal rank $n$ at points of
the intersection $S\cap \eta^{-1}(\{0\})$ that are arbitrarily
close to $0\in \C^r$. Suppose also that the ratio $(\alpha \circ
\eta:\beta \circ \eta)$ is convergent along $S$
$($in the sense of Definition {\rm\ref{convergent}}$)$. Then
$(\alpha:\beta)$ is equivalent to a nontrivial ratio of convergent
power series.
\el

\begin{proof}
Without loss of generality, $S$ is connected.
By Definition \ref{convergent}, there exist a nonnegative integer
$l$ and, for any positive integer $k$, convergent power series
$A_k(z),B_k(z)\in \C \{z\}$, $z=(z_1,\ldots,z_r)$, such that
\begin{equation}\Label{montanari}
\big(j_z^k(B_k(\alpha \circ \eta) -A_k(\beta \circ \eta))\big)|_{S}\equiv 0
\end{equation}
and $(j_z^l(A_k,B_k))|_{S}\not \equiv 0$. We may assume that
$A_k,B_k$ are convergent in a polydisc neighborhood $\Delta_k$ of
$0\in \C^r$. Choose $\nu_0 \in \N^{r}$, $|\nu_0|\leq l$, of
minimal length such that such that
\begin{equation}\Label{min-length}
(\partial^{\nu_0}_z A_l,
\partial^{\nu_0}_z B_l)|_{S}\not \equiv 0.
\end{equation}
 Then,
since $(\partial^{\nu}_z A_l,\partial^{\nu}_z B_l)|_{S}\equiv 0$ for $|\nu|<|\nu_0|$,
(\ref{montanari}) with $k=l$ implies
\begin{equation}\Label{bellissima}
\big((\partial^{\nu_0}_zB_l)(\alpha \circ
\eta)-(\partial^{\nu_0}_zA_l) (\beta\circ \eta)\big)|_S\equiv 0.
\end{equation}
By assumption on $\eta|_{S}$, we may choose a point $s_0\in S\cap \Delta_l$
arbitrarily close to $0$
with $s_0\in \Delta_l$, such that $\eta|_S$ has rank $n$ at $s_0$
and $\eta (s_0)=0$. By the rank theorem,
we may choose a right inverse of $\eta$,
$\theta\colon \Omega \to S$, holomorphic in some neighborhood $\Omega$
of $0\in \C^n$ with $\theta (0)=s_0$. Since $(\eta \circ \theta)(t)\equiv t$, we obtain from
(\ref{bellissima}) that
$((\partial^{\nu_0}_zB_l)\circ
\theta)(t) (\alpha (t))-((\partial^{\nu_0}_zA_l)\circ
\theta)(t) (\beta (t))\equiv 0$. To complete the proof of the
lemma, it remains to observe that, in view of \eqref{min-length},
$\theta$ can be chosen so that
$(((\partial^{\nu_0}_z A_l)\circ \theta) (t),
(\partial^{\nu_0}_z B_l)\circ \theta)(t))\not \equiv 0$.
\end{proof}

\subsection{Applications to pullbacks of ratios of formal power series}\Label{pullback}
The notion of convergence of a ratio of formal power series
along a submanifold introduced in Definition~\ref{convergent}
extends in an obvious way to formal power series
defined on a  complex manifold $\6X$. We also say that two ratios
$(F_1:G_1)$, $(F_2:G_2)$ of formal power series on $\6X$ are
equivalent if $F_1G_2-F_2G_1$ vanishes identically as a formal
power series on $\6X$.

Let $\6 Y$ be another complex manifold and $v\colon \6 Y\to\6 X$ a
holomorphic map defined in a neighborhood of a reference point
$y_0\in \6 Y$ with $x_0:=v(y_0)$. Consider the pullback under $v$
of a ratio $(F:G)$ of formal power series on $\6 X$ (centered at
$x_0$) and assume that it is convergent along a submanifold
$S\subset \6 Y$ through $y_0$. Under certain assumptions on the
map $v$ and on the formal power series we show in this section
that $\6 Y$ can be embedded into a larger manifold $\2{\6Y}$ and
$v$ holomorphically extended to $\2v\colon \2{\6Y}\to \6X$
such that the pullback of $(F:G)$ under $\2v$ is convergent along
a larger submanifold $\2S\subset\2{\6Y}$. The precise statement is
the following.

\bp\Label{going-up} Let $\6 X$ and $\6 Y$ be complex manifolds and
$v\colon \6 Y\to\6 X$ a holomorphic submersion with $y_0\in \6 Y$.
Let $S\subset \6 Y$ be a complex submanifold through $y_0$ and
$(F:G)$ a ratio of formal power series on $\6 X$, centered at
$x_0:=v(y_0)$, whose pullback under $v$ is convergent along $S$.
Let  $\eta \colon\6 X \to \6 C$ be a holomorphic submersion onto a
complex manifold $\6 C$. Define
$$\2{\6 Y}:=\{(y,x)\in \6 Y\times\6 X : \eta(v(y))=\eta(x)\}, \quad \2S:=\{(y,x)\in\2{\6 Y} : y\in S\},
  \quad \2v\colon \2{\6 Y}\ni (y,x)\mapsto x\in \6 X.$$
Assume that one of the following conditions hold:
\begin{enumerate}
\item[(i)] the ratio $(F:G)$ is equivalent to a
nontrivial ratio $(\alpha\circ\eta : \beta\circ\eta)$, where
$\alpha$ and $\beta$ are  formal power series on $\6 C$ centered
at $\eta(x_0)$;
\item[(ii)] the ratio $(F:G)$ is equivalent to a nontrivial ratio
of the form $(\Phi\big(Y(\eta(x)),x\big) : \Psi
\big(Y(\eta(x)),x\big))$, where $Y$ is a $\C^r$-valued formal
power series on $\6 C$ centered at $\eta(x_0)$ and $\Phi,\Psi$ are
convergent power series centered at $(Y(x_0),x_0)\in \C^r\times\6
X$.
\end{enumerate}
Then the pullback of $(F:G)$ under $\2v$ is convergent along
$\2S$. \ep

\br {\rm The conclusion of Proposition \ref{going-up} obviously
holds in the case ${\rm dim}\, \6X={\rm dim}\, \6 C$ (without
assuming neither (i) nor (ii)), and therefore, we may assume, in what
follows, that ${\rm dim}\, \6X>{\rm dim}\, \6 C$ which implies
${\rm dim}\, \2{\6Y}>{\rm dim}\, \6Y$.} \er

In order to reduce Proposition \ref{going-up} to an application of
Lemmata \ref{trivial} and \ref{nontrivial}, we need the following.

\bl\Label{keeping-asking} In the setting of Proposition {\rm
\ref{going-up}}, define $$\widehat S:=\{(y,v(y)):y\in S\}\subset
\2{\6 Y}.$$ Then the pullback of $(F:G)$ under $\2v$ is convergent
along the complex submanifold $\widehat S$. \el

The idea of the proof lies in the fact that the derivatives of
the pullbacks under $\2v$ can be expressed through
derivatives of the pullbacks under $v$ of the same power series.
For this property to hold, it is essential to assume that $v$ is submersive.

\bpf[Proof of Lemma~{\rm \ref{keeping-asking}}] The manifold $\6
Y$ can be seen as embedded into $\2{\6 Y}$ via the map $\vartheta
\colon \6 Y \ni y\mapsto (y,v(y))\in \2{\6 Y}$. Therefore, by
considering $\vartheta (\6 Y)$, we may also think of $\6 Y$ as a
submanifold in $\2{\6 Y}$. Since $v$ is a submersion, after
possibly shrinking $\6 Y$ near $y_0$ and $\6 X$ near $x_0$, we may
choose for every $y\in \6 Y$ a holomorphic right inverse of $v$,
$v_y^{-1}\colon \6 X \to \6 Y$, such that $v_y^{-1}(v(y))=y$. Such
a choice can be made by the rank theorem so that the map $\6
Y\times \6 X\ni (y,x)\mapsto v_{y}^{-1}(x)\in \6 Y$ is
holomorphic.

 Choose open neighborhoods $\Omega_1\subset \C^{{\rm dim}\, {\6Y}}$ and
 $\Omega_2\subset \C^{{\rm dim}\, \2{\6Y}-{\rm dim}\, {\6Y}}$ of the origin and local holomorphic
coordinates $(z,w)=(z(y,x),w(y,x))\in \Omega_1\times \Omega_2$ on
$\2 {\6{Y}}$ vanishing at $(y_0,x_0)\in \2{\6 Y}$ such that
$\vartheta (\6 Y)$ is given by $\{(y,x)\in \2{\6 Y}: w=0\}$.
(Hence $z|_{\6 Y}:\6 Y\ni y\mapsto z(y,v(y))\in \Omega_1$ is a
system of holomorphic coordinates for $\6Y$.) In what follows, as is
customary, we identify $S$ and $z(S)$. Since $(F:G)$ is
convergent along $S$, for any nonnegative integer $k$, there exist
convergent power series $f_k,g_k$ in $\C \{z\}$ such that
\begin{equation}\Label{have}
(j_z^{k} R_k)|_S \equiv  0, \quad R_k(z):=(F\circ v)(z) g_k(z) -
(G\circ v)(z)f_k(z),
\end{equation}
and $(j_z^l(f_k,g_k))|_S \not\equiv 0$ for some nonnegative
integer $l$ independent of $k$. In what follows we fix $k$ and may
assume, without loss of generality, that $f_k,g_k$ are holomorphic
in $\Omega_1$. We shall define convergent power series
$\2f_k,\2g_k\in \C\{z,w\}$ whose restrictions to $\{w=0\}$ are
$f_k,g_k$. For this, we set, for $z\in \Omega_1$,
$v^{-1}_z:=v^{-1}_y$ where $y\in \6 Y$ is uniquely determined by
the relation $z=z(y,v(y))$. Define holomorphic functions on
$\Omega_1\times \Omega_2$ by setting $\2f_k(z,w):=(f_k\circ
v_z^{-1}\circ \2v)(z,w)$ and $\2g_k(z,w):=(g_k\circ v_z^{-1}\circ
\2v)(z,w)$. We also set $\2R_k(z,w):=(F\circ \2v)(z,w) \2g_k(z,w)
- (G\circ \2v)(z,w)\2f_k(z,w)$. Since $v_z^{-1}$ is a right
convergent inverse for $v$, it follows from the above construction
that $\2R_k(z,w)=(R_k\circ v_z^{-1}\circ \2v)(z,w)$. Therefore, by
the chain rule, the power series mapping $j_{(z,w)}^k\2R_k$ is a
linear combination (with holomorphic coefficients in $(z,w)$) of
the components of $(j_z^kR_k)\circ (v_z^{-1}\circ \2v)$.  By
restricting to $z\in S$ and $w=0$, we obtain, in view of
(\ref{have}) and the fact that $v_z^{-1}(\2v(z,0))=z$,
\begin{equation}\Label{better}
(j_{(z,w)}^k((F\circ \2v)\2g_k-(G\circ \2v)\2f_k))|_{z\in S, \
w=0}\equiv 0.
\end{equation}
We  therefore conclude that $(\2f_k:\2g_k)$ is $k$-similar to
$(F\circ\2v:G\circ\2v)$ along $\widehat S$ since the submanifold
$\widehat S$ is given by $\{(z,0):z\in S\}$ in the
$(z,w)$-coordinates. Since $(f_k,g_k)$ is the restriction of
$(\2f_k,\2g_k)$ to $\{w=0\}$ by construction, we have
$(j_{(z,w)}^l(\2f_k,\2g_k))|_{\widehat S}\not\equiv 0$ for $l$ as
above. This shows that $(F\circ \2v:G\circ \2v)$ is convergent
along $\widehat S$ and hence completes the proof of the lemma.
\epf

\begin{proof}[Proof of Proposition {\rm \ref{going-up}}] The
statement obviously holds when $F$ and $G$ are both zero, so we
may assume that the ratio $(F:G)$ is nontrivial. Choose local
holomorphic coordinates $Z=Z(y)\in \C^{{\rm dim}\, \6 Y}$ for $\6
Y$ , vanishing at $y_0$, such that $S$ is given in these
coordinates by $\{Z=(Z^1,Z^2)\in \C^{n_1}\times \C^{n_2}:
Z^2=0\}$, with ${\rm dim}\, \6 Y=n_1+n_2$. By the construction
of $\2{\6 Y}$, we may choose holomorphic coordinates $\2 Z$ for
$\2{\6 Y}$ near $(y_0,x_0)$ of the form $\2Z=\2 Z
(y,x)=(Z(y),Z^3(y,x))\in \C^{{\rm dim}\, \6 Y}\times
\C^{n_3}$, where $Z$ is as above, $n_3={\rm dim}\, \2{\6Y}-{\rm
dim}\, \6Y$ and such that $\vartheta (\6 Y)$ is given by
$\{Z^3=0\}$. Note that the submanifolds $\widehat S$ and $\2
S$ are given in the $\2Z$-coordinates by $\{Z^2=Z^3=0\}$
and $\{Z^2=0\}$ respectively and $\eta \circ \2 v$ is
independent of $Z^3$ (again by the construction of $\2{\6Y}$).

To prove the conclusion assuming (i),
we first note that since $v$ is a submersion and
$\vartheta (\6 Y)\subset \2 {\6Y}$, it follows that $\2 v$ is a
submersion too. Therefore, the nontrivial ratio $(F\circ
\2v:G\circ \2v)$ is equivalent to the nontrivial ratio $(\alpha
\circ \eta \circ \2v:\beta \circ \eta \circ \2v)$, and this latter
is convergent along $\widehat S$ by Lemma \ref{hue} (iii) and
Lemma \ref{keeping-asking}. To complete the proof of (i), it is
enough to prove that $(\alpha \circ \eta \circ \2v:\beta \circ
\eta \circ \2v)$ is convergent along $\2 S$ (again  by Lemma
\ref{hue} (iii)). By using the $\2 Z$-coordinates for $\2{\6 Y}$
defined above, we see that the conclusion follows from a direct
application of Lemma \ref{trivial}.

The proof of the conclusion assuming (ii) follows the same lines
as above by making use of Lemma \ref{nontrivial} (instead of Lemma
\ref{trivial}). This completes the proof of Proposition
\ref{going-up}.
\end{proof}

\subsection{Pairs of submersions of finite type and meromorphic
extension}\Label{hangar}

We shall formulate our main result of this section
in terms of pairs of submersions defined on a given complex manifold.
The main example of this setting is given
by the complexification $\6M\subset\C^N\times\C^N$
of a real-analytic generic submanifold $M\subset\C^N$,
where a pair of submersions on $\6M$
is given by the projections on the first
and the last component $\C^N$ respectively.

In general, let $\6 X$, $\6 Z$ and $\6W$ be complex manifolds and $\l\colon \6
X \to \6 Z$, $\mu\colon \6 X\to \6 W$ be holomorphic submersions.
Set ${\6 X}^{(0)}:= \6 X$ and for any integer $l\ge 1$, define the
(odd) fiber product
\begin{equation}\Label{Xl}
{\6 X}^{(l)} := \{ (z_1,\ldots,z_{2l+1})\in {\6 X}^{2l+1} :
\mu(z_{2s-1})=\mu(z_{2s}),\; \l(z_{2s})=\l(z_{2s+1}),\; 1\le s\le
l \}.
\end{equation}
Analogously fiber products with even number of factors can be  defined
but will not be used in this paper.
It is easy to see that ${\6
X}^{(l)}\subset {\6 X}^{2l+1}$ is a complex submanifold. Let
\begin{equation}\Label{pi}
{\6 X}^{(l)}\ni (z_1,\ldots,z_{2l+1})\mapsto
\pi^{(l)}_j(z_1,\ldots,z_{2l+1}):=z_j\in \6 X
\end{equation} be
the restriction to ${\6 X}^{(l)}$ of the natural projection to the
$j$-th component, $1\leq j\leq 2l+1$, and denote by $\2\l\colon
{\6 X }^{(l)}\to \6 Z$ and $\2\mu\colon {\6 X}^{(l)}\to \6 W$
the maps defined by
\begin{equation}\Label{pil}
\2\l :=\l \circ \pi^{(l)}_{1},\quad \2\mu:= \mu \circ
\pi^{(l)}_{2l+1}. \end{equation}
 Then, for every $x\in\6 X$ we set
$x^{(l)}:=(x,\ldots,x)\in {\6 X}^{(l)}$ and
\begin{equation}\Label{fiber}
D_l(x):={\2\l^{-1}(\2\l(x^{(l)})}),\quad E_l(x):=
{\2\mu^{-1}(\2\mu(x^{(l)}))} \end{equation}
are complex submanifolds of ${\6X}^{(l)}$.

In the above mentioned case, i.e.\ when $\6X$
is the complexification of a real-analytic generic submanifold $M\subset\C^N$,
the construction of ${\6 X}^{(l)}$ yields the iterated complexification
$\6M^{2l}$ as defined in \cite{Z97}.
In this case the images $\2\mu(D_l(x))$ are the Segre sets
in the sense of {\sc Baouendi-Ebenfelt-Rothschild} \cite{BERacta}
and their finite type criterion says that $M$ is of finite type
in the sense of {\sc Kohn}
\cite{Kohn} and {\sc Bloom-Graham} \cite{BG}
if and only if the Segre sets of sufficiently high order
have nonempty interior. The last condition
can also be expressed in terms of ranks (see \cite{BERbook}).
Motivated by this case, we say in the above general setting
 that the pair $(\l,\mu)$ of submersions is
{\em of finite type} at a point $x_0\in\6 X$ if there exists
$l_0\geq 1$ such that the map $\2\mu_{l_0}|_{D_{l_0}(x_0)}$ has
rank equal to $\dim \6 W$ at some points of the intersection
${D_{l_0}(x_0)}\cap{E_{l_0}(x_0)}$ that are arbitrarily close to
$x_0^{(l_0)}$.

The main result of \S \ref{fou} is the following meromorphic
extension property of ratios of formal power series
that was inspired by an analogous result from \cite{MMZ1} 
in a different context. Its proof is however completely different
and will consist of repeatedly applying Proposition \ref{going-up}.

\bt\Label{mer-ext} Let $\6 X$, $\6 Z$, $\6W$ be complex manifolds
and $\l\colon \6 X \to \6 Z$, $\mu\colon \6 X\to \6 W$ be a pair
of holomorphic submersions of finite type at a point $x_0\in \6
X$. Consider formal power series $F(x),G(x)$ on $\6 X$ centered at
$x_0$ of the form $F(x)=\Phi\big(Y(\l(x)),x\big)$,
$G(x)=\Psi\big(Y(\l(x)),x\big)$, where $Y$ is a $\C^r$-valued
formal power series on $\6 W$ centered at $\l(x_0)$ and
$\Phi,\Psi$ are convergent power series on $\C^r\times\6 X$
centered at $(Y(\l (x_0)),x_0)$. Suppose that $G\not\equiv 0$ and
that  ${\6L} (F/G) \equiv 0$ holds for any holomorphic vector
field ${\6L}$ on $\6 X$ annihilating $\mu$. Then $(F:G)$ is
equivalent to a nontrivial ratio of convergent power series on $\6
X$ $($centered at $x_0)$. \et

\br\Label{kim}
 {\rm From the proof of Theorem \ref{mer-ext}, it will
follow that the ratio $(F:G)$ is even equivalent to a ratio of
the form $(\2\alpha \circ \mu:\2\beta \circ \mu)$, where
$\2\alpha,\2\beta$ are convergent power series on $\6 W$ centered
at $\mu (x_0)$ with $\2\beta \not \equiv 0$.} \er

We start by giving several lemmata
that will be used in the proof of Theorem \ref{mer-ext}.

\bl\Label{exercise}
 Let $\mu\colon \6 X\to \6 W$ be a holomorphic
submersion between complex manifolds and $F(x),G(x)$ be formal
power series on $\6 X$ centered at a point $x_0\in \6X$. Suppose
that $G\not \equiv 0$ and that $\6 L (F/G)\equiv 0$ for any
holomorphic vector field $\6 L$ on $\6 X$ that annihilates $\mu$.
Then there exist formal power series $\alpha, \beta$ on $\6 W$
centered at $\mu (x_0)$, with $\beta \not \equiv 0$, such that the
ratio $(F:G)$ is equivalent to the ratio $(\alpha \circ \mu: \beta
\circ \mu)$. \el

The proof of Lemma \ref{exercise} follows from Lemma \ref{talk}
after appropriate choices of local coordinates in $\6 X$ and
$\6W$. In the next lemma, we apply the iteration process provided by
Proposition \ref{going-up} in the context of Theorem
\ref{mer-ext}.

\bl\Label{iteration} In the setting of Theorem {\rm
\ref{mer-ext}}, the following holds. If, for some nonnegative
integer $l$, the ratio $(F\circ \pi^{(l)}_{2l+1}:G\circ
\pi^{(l)}_{2l+1})$ is convergent along $D_l(x_0)$, then the ratio
$(F\circ \pi^{(l+1)}_{2l+3}:G\circ \pi^{(l+1)}_{2l+3})$ is
convergent along $D_{l+1}(x_0)$. Here, $\pi^{(j)}_{2j+1}$ is the
projection given by {\rm (\ref{pi})} and $D_j(x_0)$ is the
submanifold given by {\rm (\ref{fiber})}, $j=l,l+1$. \el

\begin{proof}In order to apply Proposition \ref{going-up}, we first set
$\eta:=\mu$, $\6 X:=\6 X$, $\6 C:=\6 W$, $\6 Y:={\6 X}^{(l)}$,
$y_0:=x_0^{(l)}$, $v:=\pi_{2l+1}^{(l)}$ and $S:=D_l(x_0)$, where
${\6 X}^{(l)}$ and $\pi_{2l+1}^{(l)}$  are given by (\ref{Xl}) and
(\ref{pi}) respectively. Note that $v$ is a holomorphic submersion
and that, by assumption, the pullback under $v$ of the ratio
$(F:G)$ is convergent along $S$. In view of Lemma \ref{exercise},
Proposition \ref{going-up} (i) implies that, by setting
\begin{eqnarray}\nonumber
{\6 Y}_1&:=&\{(z_1,\ldots,z_{2l+1},z_{2l+2})\in {\6
X}^{(l)}\times\6 X :  \mu (z_{2l+1})=\mu (z_{2l+2}) \},\\
S_1&:=&\{(z_1,\ldots,z_{2l+1},z_{2l+2})\in {\6 Y}_1 : \l (z_1)=\l
(x_0) \},\nonumber \\ v_1&\colon &{\6 Y}_1\ni
(z_1,\ldots,z_{2l+1},z_{2l+2})\mapsto z_{2l+2}\in \6 X, \nonumber
\end{eqnarray}
the pullback of $(F:G)$ under $v_1$ is convergent along $S_1$. We
now want to apply a second time Proposition \ref{going-up}. For
this, we reset $\eta:=\l$, $\6 X:=\6 X$, $\6 C:=\6 Z$, $\6 Y:={\6
Y}_1$, $y_0:=(x_0,\ldots,x_0)\in {\6 Y}_1$, $v:=v_1$ and $S:=S_1$,
where ${\6Y}_1$, $S_1$ and $v_1$ are as above. By applying
Proposition \ref{going-up} (ii) in that context, we obtain easily
that the pullback of $(F:G)$ under $\pi^{(l+1)}_{2l+3}$ is
convergent along $D_{l+1}(x_0)$, the required conclusion.
\end{proof}

\begin{proof}[Proof of Theorem {\rm \ref{mer-ext}}]
We first claim that the pullback of $(F:G)$ under $\pi_1^{(0)}$ is
convergent along $D_0(x_0)$. Indeed, note that this is equivalent
to saying that $(F:G)$ is convergent along $\{x\in \6 X: \l (x)=\l
(x_0)\}$. By applying Proposition \ref{going-up} (ii) with
$\eta:=\l$, $\6 X:=\6 X$, $\6 C:=\6 Z$, $\6 Y:={\6X}$, $y_0:=x_0$,
$v:={\rm Id}_{\6X}$ and $S:=\{x_0\}$, and using Lemma \ref{hue}
(i), we get the desired claim. By applying Lemma \ref{iteration}
and using the finite type assumption on the pair $(\l,\mu)$, it
follows that the ratio $(F\circ \pi^{(l_0)}_{2l_0+1}:G\circ
\pi^{(l_0)}_{2l_0+1})$ is convergent along $D_{l_0}(x_0)$, where
$l_0$ is chosen so that $\2\mu|_{D_{l_0}(x_0)}$ has rank equal to
$\dim \6 W$ at some points of the intersection
${D_{l_0}(x_0)}\cap{E_{l_0}(x_0)}$ that are arbitrarily close to
$x_0^{(l_0)}$. Let $\alpha$ and $\beta$ be power series on $\6W$
given by Lemma \ref{exercise}. In view of Lemma \ref{hue} (iii),
the nontrivial ratio $(\alpha\circ \2\mu:\beta\circ \2\mu)$ is
thus convergent along $D_{l_0}(x_0)$, where $\beta \circ \2\mu\not
\equiv 0$. Since $E_{l_0}=\2\mu^{-1}(\{0\})$, one sees that Lemma
\ref{annamaria} implies that the ratio $(\alpha:\beta)$ is
equivalent to a nontrivial ratio $(\2 \alpha :\2 \beta)$ of
convergent power series on $\6 W$ (centered at $\mu (x_0)$).
Therefore, it follows from Lemma \ref{exercise} and the fact that
$\mu$ is a submersion that $(F:G)$ is equivalent to the nontrivial
ratio $(\2 \alpha \circ \mu:\2 \beta \circ \mu)$. The proof of
Theorem \ref{mer-ext} is complete.
\end{proof}

\section{Applications of Theorem~\ref{mer-ext} to ratios on generic submanifolds}\Label{realstuff}

The goal of this section is to apply the meromorphic extension
property of ratios of formal power series given by Theorem
\ref{mer-ext} to the context of real-analytic generic submanifolds
in $\C^N$, and to deduce some other properties (see Proposition
\ref{manger} below) which will be useful for the proof of the
theorems mentioned in the introduction.

 Let $M\subset \C^N$ be a real-analytic generic
submanifold of codimension $d$ through $0$, and $\rho
(Z,\1Z):=(\rho_1(Z,\1Z),\ldots,\rho_d(Z,\1Z))$ be a real-analytic
vector-valued defining function for $M$ defined in a connected
neighborhood $U$ of $0$ in $\C^N$, satisfying $\partial
\rho_1\wedge \ldots\partial \rho_d\not =0$ on $U$. Define the
complexification $\M$ of $M$ as follows
\begin{equation}\Label{comp}
\M:=\{(Z,\zeta)\in U\times U^*: \rho (Z,\zeta)=0\},
\end{equation}
where for any subset $V\subset \C^k$, we have denoted
$V^*:=\{\1w:w\in V \}$. Clearly, $\M$ is a $d$-codimensional
complex submanifold of $\C^N\times \C^N$. We say that a vector
field $X$ defined in a neighborhood of $0\in \C^N\times \C^N$ is a
$(0,1)$ vector field if it annihilates the natural projection
$\C^N\times \C^N \ni (Z,\zeta)\mapsto Z\in \C^N$. We also say that
$X$ is tangent to $\M$ if  $X(q)\in T_q\M$ for any $q\in \M$ near
the origin. We have the following consequence of
Theorem~\ref{mer-ext}.

\bt\Label{hurry} Let $M\subset \C^N$ be a real-analytic generic
submanifold through $0$ and $\M\subset \C^N_Z\times
\C^N_{\zeta}$ its complexification as given by {\rm (\ref{comp})}.
Consider formal power series $F(Z,\zeta),G(Z,\zeta)\in \C \dbl
Z,\zeta\dbr$ of the form $F(Z,\zeta)=\Phi (Y(\zeta),Z)$,
$G(Z,\zeta)=\Psi (Y(\zeta),Z)$, where $Y(\zeta)$ is a
$\C^r$-valued formal power series  and $\Phi,\Psi$ are convergent
power series centered at $(Y(0),0)\in \C^r \times \C^N$ with
$G(Z,\zeta)\not \equiv 0$ for $(Z,\zeta)\in \6M$. Suppose that $M$
is minimal at $0$ and that $\6 L (F/G)\equiv 0$ on $\6 M$ $($i.e.\
$F\6L G-G\6L F\equiv 0$ on $\6 M$$)$ for any $(0,1)$ holomorphic
vector field tangent to $\6M$. Then there exist convergent power
series $\2F(Z),\2G (Z)\in \C \{Z\}$, with $\2G(Z)\not \equiv 0$,
such that the ratios $(F:G)$ and $(\2F:\2G)$ are equivalent as
formal power series on $\6M$. \et

For the proof of the theorem, we set $\6 X:=\6 M$, $\6 Z=\6
W:=\C^N$ and define the holomorphic submersions
$$\l\colon \6 M\ni (Z,\zeta)\mapsto \zeta \in \C^N,\quad
\mu \colon \6 M\ni (Z,\zeta)\mapsto Z\in \C^N.$$

\bl\Label{julie} In the above setting, the pair $(\l,\mu)$ is of
finite type at $0\in \6M$ $($as defined in \S {\rm
\ref{hangar}}$)$ if and only if $M$ is minimal at the origin. \el

\begin{proof} For any nonnegative integer $l$, the fiber product ${\6 M}^{(l)}$ is here given by
$${\6 M}^{(l)}=\{((Z_1,\zeta_1),\ldots,(Z_{2l+1},\zeta_{2l+1}))\in {\6 M}^{2l+1}: Z_{2s-1}=Z_{2s},\ \zeta_{2s}=\zeta_{2s+1},\
1\leq s\leq l\}$$ and the maps $\2{\l}_{l}\colon {\6 M}^{(l)}\to
\C^N$, $\2\mu \colon {\6 M}^{(l)}\to \C^N$ by
$$(Z_1,\zeta_1,\ldots,Z_{2l+1},\zeta_{2l+1})\mapsto\zeta_1,\
(Z_1,\zeta_1,\ldots,Z_{2l+1},\zeta_{2l+1})\mapsto Z_{2l+1}$$
respectively.  We then have
$D_{l}(0)=\{((Z_1,\zeta_1),\ldots,(Z_{2l+1},\zeta_{2l+1}))\in {\6
M}^{(l)}:\zeta_1=0\}$ and
$E_l(0)=\{((Z_1,\zeta_1),\ldots,(Z_{2l+1},\zeta_{2l+1}))\in {\6
M}^{(l)}:Z_{2l+1}=0\}$. The reader can check that the map
$\2{\mu}_l|_{D_l(0)}$ coincides, up to a parametrization of ${\6
M}^{(l)}$, with a suitable iterated Segre mapping $v^{2l+1}$ at
$0$ as defined in \cite{BER99,BERalg}. Therefore, in view of the
minimality criterion of \cite{BERbook,BERalg} (see also \cite{BERacta}),
the pair $(\lambda,\mu)$ is of finite type at
$0\in \6M$ if and only if $M$ is minimal at $0$. The proof of the
lemma is complete.
\end{proof}

\begin{proof}[Proof of Theorem {\rm \ref{hurry}}] Since $(0,1)$ holomorphic vector fields
tangent to $\M$ coincide with holomorphic vector fields on $\M$
annihilating the submersion $\mu$, in view of Lemma \ref{julie},
we may apply Theorem \ref{mer-ext} to conclude that the ratio
$(F:G)$ is equivalent to a ratio $(\2F(Z):\2G(Z))$ of convergent
power series on $\6 M$ with $\2G(Z)\not \equiv 0$. (The fact that
$\2F,\2G$ may be chosen independent of $\zeta$ follows from Remark
\ref{kim}.) The proof is complete.
\end{proof}

In what follows, for any ring $A$, we denote, as usual, by $A[T]$,
$T=(T_1,\ldots,T_r)$, the ring of polynomials over $A$ in $r$
indeterminates. An application of Theorem \ref{hurry} is given by
the following result, which will be essential for the proof of the
main results of this paper.

\begin{Pro}\Label{manger}
Let $M\subset \C^N$ be a minimal real-analytic generic
submanifold through $0$ and $\M\subset
\C^N_{Z}\times \C^N_{\zeta}$ be its complexification given by {\rm
(\ref{comp})}. Let $F(Z):=(F_1(Z),\ldots,F_r(Z))$ be a formal
power series mapping satisfying one of the following conditions:
\begin{enumerate}
\item [(i)] there exists
$G(\zeta):=(G_1(\zeta),\ldots,G_s(\zeta))\in (\C \dbl \zeta \dbr)^s$, $G(0)=0$,
and a polynomial ${\mathcal R}(Z,\zeta,X;T)$ $\in \C\{Z,\zeta, X\}
[T]$, $X=(X_1,\ldots,X_s)$, $T=(T_1,\ldots,T_r)$, such that
${\mathcal R}(Z,\zeta,G(\zeta);T)\not \equiv 0$ for $(Z,\zeta)\in
\M$ and such that ${\mathcal R}(Z,\zeta,G(\zeta);F(Z)) \equiv  0$
for $(Z,\zeta)\in \M$;
\item [(ii)] there exists a polynomial ${\mathcal P}(Z,\zeta;\2 T,T)\in \C\{Z,\zeta\}[\2 T,T]$,
$\2 T=(\2{T}_1,\ldots,\2{T}_r)$, $T=(T_1,\ldots,T_r)$, such that
${\mathcal P}(Z,\zeta;\2 T,T)\not \equiv 0$ for $(Z,\zeta)\in \M$
and such that ${\mathcal P}(Z,\zeta;\1 F(\zeta),F(Z))\equiv 0$ for
$(Z,\zeta)\in \M$.
\end{enumerate}
Then there exists a nontrivial polynomial $\Delta (Z;T)\in
\C\{Z\}[T]$ such that $\Delta (Z,F(Z))\equiv 0$.
\end{Pro}

\begin{proof}
Let $\6R$ be as in (i) such that
\begin{equation}\Label{vanishing}
{\mathcal R}(Z,\zeta,G(\zeta);F(Z))\equiv 0,\ {\rm for}\
(Z,\zeta)\in \M.
\end{equation}
We write $\6R$ as a linear combination
\begin{equation}\Label{mon}
{\mathcal R}(Z,\zeta,G(\zeta);T)=\sum_{j=1}^{l}\delta_j(Z,\zeta,G(\zeta))\, r_{j}(T),
\end{equation}
where each $\delta_j(Z,\zeta,G(\zeta))\not \equiv 0$ for
$(Z,\zeta)\in \M$, $\delta_j(Z,\zeta,X)\in \C \{Z,\zeta\}[X]$, and
$r_j$ is a monomial in $T$. We prove the desired conclusion by
induction on the number $l$ of monomials
 in (\ref{mon}).
For $l=1$, (\ref{vanishing}), (\ref{mon}) and the fact that
$\delta_1(Z,\zeta,G(\zeta))\not \equiv 0$ on $\M$ imply that
$r_1(F(Z))\equiv 0$. Since $r_1$ is a monomial, it follows that
$F_j(Z)=0$ for some $j$ which yields the required nontrivial
polynomial identity.

Suppose now that the desired conclusion holds for any polynomial
$\6R$ whose  number of monomials is strictly less than $l$ and for
{\em any} formal power series mapping $G(\zeta)$. In view of
(\ref{vanishing}) and (\ref{mon}), we have the following identity
(understood in the field of fractions of formal power series)
\begin{equation}\Label{cms}
r_l(F(Z)) +  \sum_{j<l}
\frac{\delta_j(Z,\zeta,G(\zeta))}{\delta_l(Z,\zeta,G(\zeta))} \,
r_{j}(F(Z))\equiv 0,\ (Z,\zeta)\in \M.
\end{equation}
Let ${\mathcal L}$ be any (0,1) holomorphic vector field tangent
to $\M$. Applying ${\mathcal L}$ to (\ref{cms}) and using the fact
that ${\mathcal L}(F_j(Z))\equiv 0$ for any $j$, we obtain
\begin{equation}\Label{obtain}
\sum_{j<l}{\mathcal L}\left(
\frac{\delta_j(Z,\zeta,G(\zeta))}{\delta_l(Z,\zeta,G(\zeta))}\right)
r_{j}(F(Z)) \equiv 0,\ (Z,\zeta)\in \M.
\end{equation}
We set $Q_j(Z,\zeta):=\delta_j(Z,\zeta, G
(Z,\zeta))/\delta_l(Z,\zeta,G (Z,\zeta))$. It is easy to see that
each ratio ${\6 L} Q_j$ can be written as a ratio of the following
form
\begin{equation}\Label{form}
\frac{\2{\delta_j}(Z,\zeta,\2 G(\zeta))}{\2 {\delta_l}(Z,\zeta,\2 G(\zeta))}
\end{equation}
for some $\2 {G}(\zeta)\in (\C \dbl \zeta\dbr)^{\2 s}$ with
$\2{G}(0)=0$ and some $\2{\delta_j}(Z,\zeta,\2
X),\2{\delta_l}(Z,\zeta,\2 X)\in \C\{Z,\zeta, \2 X\}$, $\2 X\in
\C^{\2 s}$, with $\2{\delta_l}(Z,\zeta,\2 G (Z,\zeta))\not \equiv
0$ for $(Z,\zeta)\in \M$. From (\ref{obtain}), we are led to
distinguish two cases. If  for some $j\in \{1,\ldots,l-1\}$, $\6L
Q_j$ does not vanish identically on $\M$, then the required
conclusion follows from (\ref{obtain}), (\ref{form}) and the
induction hypothesis.

It remains to consider the case when ${\6L}Q_j\equiv 0$ on $\M$
for all $j$ and for all $(0,1)$ holomorphic vector fields
${\mathcal L}$ tangent to $\M$. Then each ratio $Q_j$ satisfies
the assumptions of Theorem \ref{hurry}, and therefore, there
exists $\Phi^j(Z), \Psi^j(Z)\in \C\{Z\}$ with $\Psi^j(Z)\not
\equiv 0$ such that $Q_j(Z,\zeta)=\Phi^j(Z)/\Psi^j(Z)$ for
$j=1,\ldots,l-1$. As a consequence, (\ref{cms}) can be rewritten
as
\begin{equation}\Label{cms1}
r_l(F(Z))+  \sum_{j<l} \frac{\Phi^j(Z)}{\Psi^j(Z)}\,
r_{j}(F(Z))\equiv 0.
\end{equation}
This proves the desired final conclusion  and completes the proof of
the conclusion assuming (i).

For the statement under the assumption (ii), consider a nontrivial
polynomial ${\mathcal P}(Z,\zeta;\2 {T},T)$ (on $\M$) such that
${\mathcal P}(Z,\zeta; \1{F}(\zeta),F(Z))\equiv 0$ for
$(Z,\zeta)\in \M$. We write
\begin{equation}\Label{flaw}
{\mathcal P}(Z,\zeta;\2{T},T)=\sum_{\nu \in \N^r, |\nu|\leq l}{\mathcal P}_{\nu}(Z,\zeta;\2{T}) T^{\nu},
\end{equation}
where each ${\mathcal P}_{\nu}(Z,\zeta;\2{T}) \in
\C\{Z,\zeta\}[\2{T}]$ and at least one of the $\6P_\nu$'s is
nontrivial. If there exists $\nu_0\in \N^r$ such that ${\mathcal
P}_{\nu_0}(Z,\zeta;\1{F}(\zeta))\not \equiv 0$ for $(Z,\zeta)\in
\M$, then it follows that the polynomial ${\mathcal
P}(Z,\zeta;\1{F}(\zeta),T)$ is nontrivial (on $\M$) and satisfies
${\mathcal P}(Z,\zeta;\1{F}(\zeta),F(Z))\equiv 0 $ for
$(Z,\zeta)\in \M$. Then condition (i) is fulfilled and the
required conclusion is proved above.

It remains to consider the case when ${\mathcal
P}_{\nu}(Z,\zeta;\1{F}(\zeta))\equiv 0$ on $\M$ and for any $\nu
\in \N^r$. Fix any $\nu$ such that $\6P_\nu(Z,\zeta;\2T)$ is
nontrivial for $(Z,\zeta)\in \M$.  Write $\6P_\nu(Z,\zeta;\2T)=
\sum_{|\alpha|\leq k}c_{\alpha,\nu}(Z,\zeta) \2{T}^{\alpha}$ with
each $c_{\alpha,\nu}(Z,\zeta)\in \C\{Z,\zeta\}$. Set
 $\1 {\6 P_{\nu}}(Z,\zeta;T):=\sum_{|\alpha|\leq k}\1{c_{\alpha,\nu}}(\zeta,Z)\, T^{\alpha}$. Then $\1 {\6 P_{\nu}}(Z,\zeta;T)$
 is a nontrivial polynomial (on $\M$) and satisfies $\1{\6 P_{\nu}}(Z,\zeta;F(Z))\equiv 0$ for $(Z,\zeta)\in \M$.
 Here again, condition (i) is fulfilled
and the desired conclusion follows. The proof is complete.
\end{proof}

We conclude by mentioning the following result proved in
\cite[Theorem 5.1]{M01} and which is an immediate consequence of
Proposition \ref{manger} (i) and Proposition \ref{recall}.

\begin{Cor}\Label{note}
Let $M\subset \C^N$ be a  minimal real-analytic generic
submanifold through the origin, $\M\subset \C^N_{Z}\times
\C^N_{\zeta}$ its complexification as given by {\rm (\ref{comp})}
and $F(Z)\in \C\dbl Z\dbr$. Assume that there exists
$G(\zeta):=(G_1(\zeta),\ldots,G_s(\zeta))\in (\C \dbl \zeta
\dbr)^s$ with $G(0)=0$ and a polynomial ${\mathcal
R}(Z,\zeta,X;T)\in \C\{Z,\zeta,X\}[T]$, $X=(X_1,\ldots,X_s)$,
$T\in \C$, such ${\mathcal R}(Z,\zeta,G(\zeta);T)\not \equiv 0$
for $(Z,\zeta)\in \M$ and such that ${\mathcal
R}(Z,\zeta,G(\zeta);F(Z))\equiv 0$ for $(Z,\zeta)\in \M$. Then
$F(Z)$ is convergent.
\end{Cor}

\section{Zariski closure of the graph of a formal map}\Label{Zar}
Throughout this section, we let $f\colon (\C_Z^N,0)\to
(\C_{Z'}^{N'},0)$ be a formal map. As in \S \ref{formalstuff}, we
associate to $f$ its graph $\Gamma_f\subset \C_{Z}^N\times
\C_{Z'}^{N'}$ seen as a formal submanifold.  Given a (germ at
$(0,0)\in \C^N\times \C^{N'}$ of a) holomorphic function
$H(Z,Z')$, we say that $H$ vanishes on $\Gamma_f$ if the formal
power series $H(Z,f(Z))$ vanishes identically. If $A\subset
\C^N\times \C^{N'}$ is a (germ through the origin of a)
complex-analytic subset, we further say that the graph of $f$ is
contained in $A$, and write $\Gamma_f\subset A$, if any (germ at
$(0,0)\in \C^N\times \C^{N'}$ of a) holomorphic function $H(Z,Z')$
that vanishes on $A$, vanishes also on $\Gamma_f$. The goal of
this section is to define and give some basic properties of the
Zariski closure of the graph $\Gamma_f\subset \C^N\times \C^{N}$
over the ring $\C\{Z\}[Z']$.

\subsection{Definition}\Label{defzar}
For $f$ as above, define the {\em  Zariski closure} of $\Gamma_f$
with respect to the ring $\C\{Z\}[Z']$ as the germ ${\6
Z}_f\subset \C^N\times \C^{N'}$ at $(0,0)$ of a complex-analytic
set defined as the zero-set  of all elements in $\C \{Z\}[Z']$
vanishing on $\Gamma_f$. Note that since $\6Z_f$ contains the
graph of $f$, it follows that ${\rm dim}_{\C}\, \6Z_f \geq N$. In
what follows, we shall denote by $\mu (f)$ the dimension of the
Zariski closure $\6Z_f$. Observe also that since the ring $\C \dbl
Z\dbr$ is an integral domain, it follows that $\6Z_f$ is
irreducible over $\C\{Z\}[Z']$.

\subsection{Link with transcendence degree}\Label{link}
In this section, we briefly discuss a link between the dimension
of the Zariski closure $\mu (f)$ defined above and the
 transcendence degree of a certain field extension.
The reader is referred to \cite{ZS} for basic notions from field
theory used here. In what follows, if $\K\subset \LL$ is a field
extension and $(x_1,\ldots,x_l)\in (\LL)^l$, we write $\K
(x_1,\ldots,x_l)$ for the subfield of $\LL$
 generated by $\K$ and $(x_1,\ldots,x_l)$.

We denote by $\MM_N$ the quotient field  of the ring $\C\{Z\}$ and
consider the field extension ${\MM}_N\subset
{\MM}_N(f_1(Z),\ldots,f_{N'}(Z))$ where we write
$f(Z)=(f_1(Z),\ldots,f_{N'}(Z))$. We then define the {\em
transcendence degree} of the formal map $f$, denoted in what
follows by $m(f)$, to be the transcendence degree of the above
finitely generated field extension. (We should point out that this
notion of transcendence degree of a formal map is in general
different from the one discussed in \cite{M00, M01}.) We have the
following standard relation between $m (f)$ and $\mu (f)$:

\begin{Lem}\Label{hp}
For any formal map $f\colon (\C^N,0)\to (\C^{N'},0)$, one has $\mu (f)=N+m(f)$.
\end{Lem}

The following well-known proposition shows the relevance of $\mu (f)$
 for the study of the convergence of the map $f$.

\begin{Pro}\Label{bof}
Let $f\colon (\C^N,0)\to (\C^{N'},0)$ be a formal map and $\mu (f)$ as above.
Then the following are equivalent:
\begin{enumerate}
\item [(i)] $\mu (f)=N$;
\item [(ii)] $f$ is convergent.
\end{enumerate}
\end{Pro}

Proposition \ref{bof} is a straightforward consequence of
Proposition \ref{recall}.

\section{Local geometry of the Zariski closure}\Label{NEW1}

In this section, we keep the notation of \S \ref{Zar}.
Our goal is to study the Zariski closure defined in \S
\ref{defzar} near some points of smoothness. It is worth
mentioning the striking analogy of the approach used here with
that of \cite{Pu1,Pu2,CPS00,MMZ01} for the study of the analytic regularity of
${\mathcal C}^{\infty}$-smooth CR-mappings.

\subsection{Preliminaries}\Label{momo}
Throughout \S \ref{NEW1}, we assume that the dimension of the
Zariski closure $\6Z_f$  satisfies
\begin{equation}\Label{golf}
\mu(f)<N+N'.
\end{equation}
In what follows, for an open subset $\Omega \subset \C^k$, we
denote by $\6O (\Omega)$ the ring of holomorphic functions in
$\Omega$. Recall also that we use the notation $\Omega^*$ for the
subset $\{\1{q}:q\in \Omega\}$.

In \S \ref{Zar}, we saw that $\mu (f)\ge N$ and $m:=m (f)=\mu
(f)-N$ coincides with the transcendence degree of the field
extension ${\MM}_{N} \subset {\MM}_N (f_1(Z),\ldots,f_{N'}(Z))$,
where $f(Z)=(f_1(Z),\ldots,f_{N'}(Z))$. As a consequence, there
exist integers $1\leq j_1<\ldots<j_m< N'$ such that
$f_{j_1}(Z),\ldots,f_{j_m}(Z)$ form a transcendence basis of
${\MM}_N(f_1(Z),\ldots,f_{N'}(Z))$ over ${\MM}_N$. After
renumbering the coordinates $Z':=(z',w')\in \C^m\times \C^{N'-m}$
and setting $m':=N'-m$, we may assume that
\begin{equation}\Label{heure}
 f=(g,h)\in \C_{z'}^m\times \C_{w'}^{m'},
\end{equation}
where $g=(g_1,\ldots,g_m)$ forms a transcendence basis of
$\MM_N(f_1,\ldots,f_{N'})$ over ${\MM}_N$.

Since the components of the formal map $h:(\C_Z^N,0)\to
(\C_{w'}^{m'},0)$ are algebraically dependent over ${\MM}_{N}(g)$,
there exist
 monic polynomials
$P_j(T)\in {\MM}_N(g)[T]$, $j=1,\ldots,m'$, such that if $h=(h_1,\ldots,h_{m'})$, then
\begin{equation}\Label{base}
P_j(h_j)=0,\ j=1,\ldots,m',\ {\rm in}\ {\MM}_N(f).
\end{equation}
As a consequence, there exist non-trivial polynomials $\widehat {P_j}(T)\in \C\{Z\}[g][T]$,
$j=1,\ldots,m'$, such that
\begin{equation}\Label{ase}
\widehat{P_j}(h_j)=0,\ j=1,\ldots,m'.
\end{equation}
For every $j=1,\ldots,m'$, we can write
\begin{equation}\Label{display1}
\widehat{P_j}(T)=\sum_{\nu \leq k_j}q_{j\nu}T^{\nu},
\end{equation}
where each $q_{j\nu}\in \C \{Z\}[g]$, $q_{jk_j}\not\equiv 0$ and $k_j\geq 1$. Since
each $q_{j\nu}$ is in $\C\{Z\}[g]$, we can also write
\begin{equation}\Label{display2}
q_{j\nu}=q_{j\nu}(Z)=R_{j\nu}(Z,g(Z))
\end{equation}
where $R_{j\nu}(z,z')\in \C\{Z\}[z']$.

Let $\Delta_0^N$ be a polydisc neighborhood of $0$ in $\C^N$ such that
 the Zariski closure $\6Z_f$
can be represented by an irreducible (over the ring $\C\{Z\}[Z']$)
closed analytic subset of $\Delta_0^N\times\C^{N'}$ (also denoted
by $\6Z_f$). We have the inclusion
\begin{equation}\Label{inclusion}
\Gamma_f\subset \6Z_f\subset \C^N\times\C^{N'}.
\end{equation}
Define
\begin{equation}\Label{write}
\2{P_j}(Z,z';T):=\sum_{\nu=0}^{k_j}R_{j\nu}(Z,z')T^{\nu}\in {\mathcal O}(\Delta_{0}^N)[z'][T],\ j=1,\ldots,m'.
\end{equation}
It follows from (\ref{ase}) -- (\ref{display2}) that one has
\begin{equation}\Label{change}
\2{P_j}(Z,g(Z);h_j(Z))\equiv 0,\ {\rm in}\ \C\dbl Z\dbr, \
j=1,\ldots,m'.
\end{equation}
Here each $R_{j\nu}(Z,z')\in {\mathcal O}(\Delta_0^N)[z']$, $k_j\geq 1$, and
\begin{equation}\Label{1}
R_{jk_j}(Z,g(Z))\not \equiv 0.
\end{equation}
Moreover, since $\C \{Z\}[z'][T]$ is a unique factorization domain
(see e.g.\ \cite{ZS}), we may assume that the polynomials given by (\ref{write}) are irreducible.

Consider the complex-analytic variety $\6V_f\subset \C^N\times \C^{N'}$
through $(0,0)$ defined by
\begin{equation}\Label{subvariety}
\6V_f:=\{(Z,z',w')\in \Delta_0^N\times \C^{m}\times \C^{m'}:\2{P_j}(Z,z';w'_j)=0,\ j=1,\ldots,m'\}.
\end{equation}
By (\ref{change}), $\6V_f$ contains the graph $\Gamma_f$ and hence
the Zariski closure $\6Z_f$. In fact, since by Lemma \ref{hp},
$\dim_{\C} \6Z_f=\mu_p(f)=N+m$, it follows from the construction
that $\6Z_f$ is the (unique) irreducible component of $\6V_f$
(over $\C\{Z\}[Z']$) containing $\Gamma_f$. Note that $\6V_f$ is
not irreducible in general and, moreover, may have a dimension
larger than $\mu (f)$.

For $j=1,\ldots,m'$, let $\2{D_j}(Z,z')\in {\mathcal O}(\Delta_0^N)[z']$ be the
discriminant of the polynomial $\2{P_j}(Z,z';T)$ (with respect to $T$). Consider the complex-analytic set
\begin{equation}\Label{D}
\2{\mathcal D}:= \cup_{j=1}^{m'}\{(Z,z')\in \Delta_0^N \times \C^{m}: \2{D_j}(Z,z')=0\}.
\end{equation}
 By the irreducibility of each polynomial $\2{P_j}(Z,z';T)$, we
 have $\2{D_j}(Z,z')\not \equiv 0$ in $\Delta_0^N\times \C^m$,
 for $j=1,\ldots,m'$. Therefore from the algebraic independence of the components of the formal map $g$
over ${\MM}_N$, it follows that the graph of $g$ is not (formally) contained in $\2{\mathcal D}$, i.e.\ that
\begin{equation}\Label{discrim}
\2{D_j}(Z,g(Z))\not \equiv 0,\ {\rm for}\ j=1,\ldots,m'.
\end{equation}
We also set
\begin{equation}\Label{sigma}
\6E:=\cup_{j=1}^{m'}\{(Z,z')\in \Delta_0^N\times \C^m: R_{jk_j}(Z,z')=0\}.
\end{equation}
It is well-known that $\6E \subset \2{\6D}$, and hence the graph of $g$ is not contained in $\6E$ too.

\subsection{Description  near smooth points}\Label{curves}
By the implicit function theorem, for any point $(Z_0,Z_0')\in
\6V_f$, $Z_0'=(z_0',w_0')\in \C^m\times \C^{m'}$, with
$(Z_0,z_0')\not \in \2{\6{D}}$, there exist polydisc neighborhoods
of $Z_0$, $z_0'$ and $w_0'$, denoted by $\Delta_{Z_0}^N \subset
\Delta_0^N\subset \C^N$, $\Delta^m_{z_0'}\subset \C^{m}$,
$\Delta_{w_0'}^{m'}\subset \C^{m'}$ respectively and a holomorphic
map
\begin{equation}\Label{vue}
\theta(Z_0,Z_0';\cdot):\Delta_{Z_0}^N \times \Delta_{z_0'}^m\to \Delta_{w_0'}^{m'}
\end{equation} such that
for $(Z,z',w')\in  \Delta_{Z_0}^N\times
\Delta_{z_0'}^m\times \Delta_{w_0'}^{m'}$,
\begin{equation}\Label{caract}
(Z,z',w')\in \6V_f  \iff  w'=\theta(Z_0,Z_0';Z,z').
\end{equation}
Note that if moreover $(Z_0,Z_0')\in \6Z_f$, then (\ref{caract}) is equivalent to  $(Z,z',w')\in \6Z_f $.
For any point $(Z_0,Z_0')\in  \6Z_f$ with $(Z_0,z_0')\not \in \2{\6{D}}$, consider the complex submanifold $\6Z_f(Z_0,Z_0')$
defined by setting
\begin{equation}\Label{coiffeur}
\6Z_f(Z_0,Z_0'):=\6Z_f\cap (\Delta_{Z_0}^N\times \Delta_{z_0'}^m\times \Delta_{w_0'}^{m'}).
\end{equation}
Note that for any point $(Z_0,Z_0')$ as above,
  by making the holomorphic change of coordinates $(\2 Z,\2 Z')=(Z,\phi (Z,Z'))\in \C^N\times \C^{N'}$ where
$\phi (Z,Z')=\phi (Z,(z',w')):=(z',w'-\theta (Z_0,Z_0';Z,z'))$,
the submanifold $\6Z_f(Z_0,Z_0')$ is given  in these new coordinates by
\begin{equation}\Label{complex}
\6Z_f(Z_0,Z_0')=\{(\2 Z,\2 Z')\in \Delta_{Z_0}^N\times \Delta_{z_0'}^m\times \C^{m'}: \2 Z_{m+1}'=\ldots=\2 Z_{N'}'=0\},
\end{equation}
where we write $\2 Z'=(\2 Z'_1,\ldots,\2 Z'_{N'})$.

We summarize the above in the following proposition.

\begin{Pro}\Label{summa}
Let $f\colon (\C^N,0)\to (\C^{N'},0)$ be a formal map and $\6Z_f$
the Zariski closure of the graph of $f$ as defined in {\rm \S
\ref{defzar}}. Suppose that $\mu (f)<N+N'$. Then for any point
$(Z_0,Z_0')\in \6Z_f$ with $(Z_0,z_0')\not \in \2{\6D}$, where $\2
{\6D}$  is given by {\rm (\ref{D})}, there exists a holomorphic
change of coordinates near $(Z_0,Z_0')$ of the form $(\2 Z,\2
Z')=(Z,\phi (Z,Z'))\in \C^N\times \C^{N'}$ such that the complex
submanifold  $\6Z_f$ is given near $(Z_0,Z_0')$ by {\rm
(\ref{complex})}, with $m=\mu (f)-N$.
\end{Pro}

For $(Z_0,Z_0')\in \6V_f$ with $(Z_0,z_0')\not \in \2{\6D}$ and
$(\zeta,\chi')\in (\Delta_{Z_0}^N)^*\times (\Delta_{z_0'}^m)^*$,
we define the $\C^{m'}$-valued holomorphic map
\begin{equation}\Label{thetabar}
\1{\theta}(Z_0,Z_0';\zeta,\chi'):= \1{\theta(Z_0,Z_0';\1
\zeta,\1{\chi'})},
\end{equation}
where $\theta(Z_0,Z_0';\cdot)$ is given by (\ref{vue}). The
following lemma will be important for the proof of Theorem
\ref{straight} below.

\begin{Lem}\Label{hope}
With the above notation, for any polynomial $r(Z',\zeta')\in \C
[Z',\zeta']$, $(Z',\zeta')\in \C^{N'} \times \C^{N'}$, there
exists a  nontrivial polynomial $\6R_0(Z,\zeta,z',\chi';T)\in
{\mathcal O}(\Delta_{0}^N\times \Delta_{0}^N)[z',\chi'][T]$, $T\in
\C$, such that for any point $(Z_0,Z_0')\in \6V_f$ with
$(Z_0,z_0')\not \in \2{\6D}$, one has
\begin{equation}\Label{dima}
\6R_0\big(Z,\zeta,z',\chi';r\big(z',\theta(Z_0,Z_0';Z,z'),
\chi',\1{\theta}(Z_0,Z_0';\zeta,\chi')\big)\big)\equiv 0,
\end{equation}
for $(Z,z')\in \Delta_{Z_0}^N\times \Delta_{z_0'}^m$ and
$(\zeta,\chi')\in (\Delta_{Z_0}^N)^*\times (\Delta_{z_0'}^m)^*$.
Moreover, $\6R_0$ can be chosen with the following property: for
any real-analytic generic submanifold $M\subset \C^N$ through the
origin, $\6R_0(Z,\zeta,z',\chi';T)\not \equiv 0$ for
$((Z,\zeta),z',\chi',T)\in (\6M\cap (\Delta_0^N\times
\Delta_0^N))\times \C^m\times \C^m\times \C$ , where $\M$ is the
complexification of $M$ as defined by {\rm (\ref{comp})}.
\end{Lem}

\begin{proof}
For $(Z,z')\in \Delta_0^N\times \C^m$ with $(Z,z')\not \in \6E$, where $\6E$ is given by (\ref{sigma}), and for
$j=1,\ldots,m'$,
 we denote by $\sigma^{(j)}_1(Z,z'),\ldots,\sigma^{(j)}_{k_j}(Z,z')$ the $k_j$ roots (counted with multiplicity) of
 the polynomial $\2{P_j}(Z,z';T)$ given by
 (\ref{write}). Similarly,  for $(\zeta,\chi')\in \Delta_0^N\times \C^m$
 with $(\zeta,\chi')\not \in \6E$, $\1 \sigma^{(j)}_1(\zeta,\chi'),\ldots,\1 \sigma^{(j)}_{k_j}(\zeta,\chi')$ denote the $k_j$ roots of the
 polynomial $\1{\2{P_j}}(\zeta,\chi';T):=\sum_{\nu=0}^{k_j}\1 {R_{j\nu}}(\zeta,\chi')T^{\nu}$ (obtained from (\ref{write})).
 (Note for any $j=1,\ldots, m'$ and for any $1\leq \nu \leq k_j$ one has $\1 \sigma^{(j)}_{\nu}(\zeta,\chi')=
 \1{\sigma^{(j)}_{\nu}(\1{\zeta},\1{\chi'})}$, which justifies the slight abuse of notation made here.)
 Fix $r(Z',\zeta')\in \C [Z',\zeta']$ and set for $(Z,z')$ and
$(\zeta,\chi')$ as above
\begin{multline}\Label{trick}
\6R_1(Z,\zeta,z',\chi';T):=\\
\prod_{l_1=1}^{k_1}\ldots \prod_{l_{m'}=1}^{k_{m'}} \prod_{n_1=1}^{k_1}\ldots \prod_{n_{m'}=1}^{k_{m'}}
\Big(T- r\big(z',\sigma^{(1)}_{n_1}(Z,z'),\ldots,\sigma^{(m')}_{n_{m'}}(Z,z'),
\chi',\1 \sigma^{(1)}_{l_1}(\zeta,\chi'),\ldots, \1 \sigma^{(m')}_{l_{m'}}(\zeta,\chi')\big)\Big).
\end{multline}
It follows from Newton's theorem that (\ref{trick}) may be rewritten as
\begin{equation}\Label{risk}
\6R_1(Z,\zeta,z',\chi';T)=T^{\delta}+\sum_{\nu <\delta}A_{\nu}(Z,\zeta,z',\chi')T^{\nu},
\end{equation}
for some positive integer $\delta$, and where $A_{\nu}$ is of the form
\begin{equation}\Label{eps}
A_{\nu}(Z,\zeta,z',\chi')= B_{\nu}
\bigg(z',\chi', \bigg( \Big(\frac{R_{j\alpha}(Z,z')}{R_{jk_j}(Z,z')}\Big)_{\alpha \leq k_j}
,\Big(\frac{\1{R_{j\beta}}(\zeta,\chi')}{\1{R_{jk_j}}(\zeta,\chi')}\Big)_{\beta \leq k_j}\bigg)_{0\leq j\leq m'}
\bigg),
\end{equation}
with $B_{\nu}$ being polynomials in their arguments (depending only
on the coefficients of $r(Z',\zeta')$). In view of (\ref{risk}) and (\ref{eps}), it is
clear that there exists $C(Z,\zeta,z',\chi')\in {\mathcal
O}(\Delta_{0}^N\times \Delta_{0}^N)[z',\chi']$ with
$C(Z,\zeta,z',\chi')\not \equiv 0$ such that
\begin{equation}\Label{fr}
\6R_0(Z,\zeta,z',\chi';T):=C(Z,\zeta,z',\chi')\cdot \6R_1(Z,\zeta,z',\chi';T)
\in {\mathcal O}(\Delta_{0}^N\times \Delta_{0}^N)[z',\chi'][T].
\end{equation}
($C$ is obtained by clearing denominators in (\ref{eps}) for all
$\nu<\delta$, and hence is a product of two nonzero terms, one in
the ring ${\mathcal O}(\Delta_{0}^N)[z']$ and the other in
${\mathcal O}(\Delta_{0}^N)[\chi']$.) Since for any fixed
$(Z_0,Z_0')\in \6V_f$ with $(Z_0,z_0')\not \in \2{\6D}$ and for
any $(Z,z')\in \Delta_{Z_0}^N\times \Delta_{z_0'}^m$, $(Z,z')\not
\in \6E$ and the $j$-th component of $\theta (Z_0,Z_0'; Z,z')$ is
a root of the polynomial $\2{P_j}(Z,z';T)$ by (\ref{caract}) and
(\ref{subvariety}), it follows that
 $\6R_0$ satisfies (\ref{dima}). Finally, the last desired property of $\6R_0$ is easily seen from the explicit
 construction of the polynomial,
i.e.\ from the fact that $C(Z,\zeta,z',\chi')$ cannot vanish
identically when restricted to $(\6M\cap (\Delta_0^N\times
\Delta_0^N))\times \C^m\times \C^m$. The proof of Lemma \ref{hope}
is complete.
\end{proof}

\subsection{Approximation by convergent maps}\Label{beret}
Since the graph of the formal map $f$ is contained in ${\mathcal
Z}_f$, by applying {\sc Artin}'s approximation theorem \cite{A68},
for any nonnegative integer $\kappa$, there exists a convergent
map $f^{\kappa}:(\C^N,0)\to (\C^{N'},0)$ agreeing with $f$ at $0$
up to order $\kappa$ such that the graph of $f^{\kappa}$ is
contained in $\6Z_f$. We may assume that the  maps $f^{\kappa}$
are convergent in a polydisc neighborhood
$\Delta^N_{0,\kappa}\subset \Delta^N_{0}$ of $0$ in $\C^N$.
Following the splitting (\ref{heure}), we write
$f^{\kappa}=(g^{\kappa},h^{\kappa})$ and set
\begin{equation}\Label{hungry}
\Sigma^{\kappa}:= \{Z\in \Delta^N_{0,\kappa}:
(Z,g^{\kappa}(Z))\not \in \2{\6D} \},\quad
\Gamma_{f^{\kappa}}:=\{(Z,f^{\kappa}(Z)):Z\in \Delta^N_{0,\kappa}\}\subset
\C^N\times \C^{N'}.
\end{equation}
Observe that since $\Gamma_g$ is not contained in $\2{\6D}$ (see
\S \ref{momo}), it follows that for $\kappa$ large enough, say
$\kappa\geq \2{\kappa}$, the graph of $g^{\kappa}$ is not
contained in $\2{\6D}$ too, and therefore
$\Delta^N_{0,\kappa}\setminus \Sigma^{\kappa}$ is dense in
$\Delta^N_{0,\kappa}$. We may therefore, in what follows, assume
that $\2{\kappa}=0$. Note also that since the graph of
$f^{\kappa}$ is contained in $\6Z_f$, in view of (\ref{caract}),
one has for any $Z_0\in \Delta^N_{0,\kappa}\setminus
\Sigma^{\kappa}$
\begin{equation}\Label{mmz}
h^{\kappa}(Z)=\theta (Z_0,f^{\kappa}(Z_0);Z,g^{\kappa}(Z)),
\end{equation}
for all $Z$ in some (connected) neighborhood
$\Omega^{\kappa}_{Z_0}\subset \Delta_{Z_0}^N\cap \Delta_{0,\kappa}^N$ of $Z_0$.

\section{Main technical result}\Label{embed}

With all the tools defined in \S \S \ref{Zar}--\ref{NEW1} at our
disposal, we are now ready to prove the following statement from
which all theorems mentioned in the introduction will follow. In
what follows, we keep the notation introduced in \S \ref{Zar}--\S
\ref{NEW1}.

\begin{Thm}\Label{straight}
Let $f\colon (\C^N,0)\to (\C^{N'},0)$ be a formal map, $\6Z_f$ the
Zariski closure of $\Gamma_f$ as defined in {\rm \S \ref{defzar}}
and $(f^{\kappa})_{\kappa \geq 0}$
 the convergent maps given in {\rm \S \ref{beret}} $($associated to $f$ and ${\6Z}_f$$)$. Let $M\subset \C^N$ be a minimal
real-analytic generic submanifold  through the origin. Assume that
$f$ sends $M$ into $M'$ where $M'\subset \C^{N'}$ is a proper
real-algebraic subset through the origin. Then, shrinking $M$
around the origin if necessary, there exist a positive integer
${{\kappa}}_0$ and an appropriate union $Z_f$ of local real-analytic
irreducible components of ${\6Z}_f\cap (M\times \C^{N'})$
such that the following hold:
\begin{enumerate}

\item[(i)] $\mu (f)<N+N'$ for $\mu (f)={\rm dim}\, {\6Z}_f$;
\item[(ii)] for any
$\kappa\geq {\kappa}_0$, $\Gamma_{f^{\kappa}}\cap (M\times \C^{N'})\subset Z_f\subset M\times M'$,
where $\Gamma_{f^{\kappa}}$ is given by {\rm (\ref{hungry})};
\item[(iii)] $Z_f$ satisfies the following straightening
property: for any $\kappa \geq \kappa_0$, there exists a
neighborhood $M^{\kappa}$ of $0$ in $M$ such that for any point
$Z_0$ in a dense open subset of $M^{\kappa}$, there exists a
neighborhood $U_{Z_0}^{\kappa}$ of $(Z_0,f^{\kappa}(Z_0))$ in
$\C^N\times \C^{N'}$ and a holomorphic change of coordinates in
$U_{Z_0}^{\kappa}$ of the form $(\2Z,\2Z')=\Phi^{\kappa}
(Z,Z')=(Z,\phi^{\kappa}(Z,Z'))\in \C^N\times \C^{N'}$ such that
\begin{equation}\Label{bisreal}
Z_f\cap U^{\kappa}_{Z_0} = \{(Z,Z')\in U^{\kappa}_{Z_0} : Z\in
M, \; \2Z'_{m+1}=\cdots=\2Z'_{N'}=0\},
\end{equation}
where $m=\mu(f)-N$.
\end{enumerate}
\end{Thm}

For the proof of the above result, we shall need the following key
proposition.

\begin{Pro}\Label{bourde}
Under the assumptions of Theorem {\rm \ref{straight}}, shrinking
$M$ around the origin if necessary, the following hold:
\begin{enumerate}
\item[(i)] $\mu (f)<N+N'$;
\item[(ii)] there exists a positive integer $\kappa_0$ such that for all $\kappa \geq \kappa_0$ and for all points
$Z_0\in (M\cap \Delta^{N}_{0,\kappa})\setminus \Sigma^{\kappa}$,
the real-analytic submanifold $\6Z_f((Z_0,f^{\kappa}(Z_0))\cap
(M\times \C^{N'})$ is contained in $M\times M'$. Here
$\Sigma^{\kappa}$ and $\6Z_f((Z_0,f^{\kappa}(Z_0)))$  are given by
{\rm (\ref{hungry})} and {\rm (\ref{coiffeur})} respectively and
$\Delta^{N}_{0,\kappa}$ is a polydisc of convergence of
$f^{\kappa}$.
\end{enumerate}
\end{Pro}

\begin{proof}[Proof of Proposition {\rm \ref{bourde} (i)}]
Since $M'$ is a proper real-algebraic subset of $\C^{N'}$, there
exists a nontrivial polynomial $\rho'(Z',\1{Z'})\in \C[Z',\1{Z'}]$
vanishing on $M'$. By assumption, $f$ sends $M$ into $M'$ and
therefore we have $\rho'(f(Z),\overline{f(Z)})\equiv 0$ for $Z\in
M$, or, equivalently
\begin{equation}\Label{goal}
\rho'(f(Z),\1{f}(\zeta))=0,\ (Z,\zeta)\in \M,
\end{equation}
where $\M$ is the complexification of $M$ as given by
(\ref{comp}). It follows from Proposition \ref{manger} (ii)
(applied to $F(Z):=f(Z)=(f_1(Z),\ldots,f_{N'}(Z))$) that the
components $f_1(Z),\ldots,f_{N'}(Z)$ satisfy a nontrivial
polynomial identity with coefficients in $\C \{Z\}$. This implies
that $\mu (f)<N+N'$. The proof of Proposition \ref{bourde} (i) is
complete.
\end{proof}

By Proposition \ref{bourde} (i), we may now assume that \eqref{golf}
holds and hence the arguments of \S \ref{NEW1} apply. Since $M'$ is a
real-algebraic subset of $\C^{N'}$, it is given by
\begin{equation}\Label{orleans}
 M':=\{Z'\in \C^{N'}:\rho'_1(Z',\1{Z'})=\ldots=\rho'_{l}(Z',\1{Z'})=0\},
 \end{equation}
 where each $\rho'_j(Z',\1{Z'})$, for $j=1,\ldots,l$, is a real-valued polynomial in $\C[Z',\1{Z'}]$.

\begin{proof}[Proof of Proposition {\rm \ref{bourde} (ii)}.]
By shrinking $M$ around the origin, we may assume that
$M$ is connected and is contained in $\Delta_{0}^N$. We proceed by
contradiction. Then, in view of (\ref{caract}), (\ref{coiffeur})
and (\ref{orleans}), there exists $j_0\in \{1,\ldots,l\}$ and a
subsequence $(f^{s_k})_{k\geq 0}$ of $(f^{\kappa})_{\kappa\geq 0}$
such that for any $k$, there exists $Z^{k}\in M\cap
\Delta^N_{0,s_k}$ such that
\begin{equation}\Label{colloq}
\rho'_{j_0}\big(z',\theta \big(Z^{k},f^{s_k}(Z^k);Z,z'\big),
\1{z'},\1{\theta\big(Z^{k},f^{s_k}(Z^k);Z,z'\big)}\big)\not \equiv 0,
\end{equation}
for $(Z,z')\in (M\cap \Delta_{Z^k}^N)\times
\Delta_{g^{s_k}(Z^k)}$. After complexification of (\ref{colloq}),
we obtain
\begin{equation}\Label{col}
\rho'_{j_0}\big(z',\theta\big(Z^{k},f^{s_k}(Z^k);Z,z'\big),
\chi',\1{\theta}\big(Z^{k},f^{s_k}(Z^k);\zeta,\chi'\big)\big)\not \equiv 0,
\end{equation}
for $(Z,\zeta)\in \M\cap (\Delta_{Z^k}^N\times (\Delta_{Z^k}^N)^*)$ and
$(z',\chi')\in \Delta_{g^{s_k}(Z^k)}\times (\Delta_{g^{s_k}(Z^k)})^*$. By Lemma \ref{hope} applied to
$r(Z',\zeta'):=\rho'_{j_0}(Z',\zeta')$, there exists a nontrivial polynomial
$\6R_0(Z,\zeta,z',\chi';T)\in {\mathcal O}(\Delta_{0}^N\times \Delta_{0}^N)[z',\chi'][T]$ such that for any positive integer $k$ one has
\begin{equation}\Label{trio}
\6R_0\big(Z,\zeta,z',\chi';\rho'_{j_0}\big(z',\theta\big(Z^k,f^{s_k}(Z^k); Z,z'\big),
\chi',\1{\theta}\big(Z^k,f^{s_k}(Z^k);\zeta,\chi'\big)\big)\big)\equiv 0,
\end{equation}
for $(Z,\zeta)$ and $(z',\chi')$ as above. By Lemma \ref{hope},
$\6R_0$ does not vanish identically when restricted to $(\6M\cap
(\Delta_0^N\times \Delta_0^N))\times \C^m\times \C^m \times \C$
and therefore we may write
\begin{equation}\Label{blabla}
\6R_0(Z,\zeta,z',\chi';T)=T^{\eta}\cdot
\6R_{00}(Z,\zeta,z',\chi';T),
\end{equation}
for $((Z,\zeta),z',\chi',T)$ as above and for some integer $\eta$
and some polynomial $\6R_{00}(Z,\zeta,z',\chi';T)\in {\mathcal
O}(\Delta_{0}^N\times \Delta_{0}^N)[z',\chi'][T]$ satisfying
\begin{equation}\Label{nonde}
\6R_{00}(Z,\zeta,z',\chi';0)\not \equiv 0,\
((Z,\zeta),z',\chi')\in \M \times \C^m\times \C^m.
\end{equation}
We also write
\begin{equation}\Label{ji}
\6R_{00}(Z,\zeta,z',\chi';T)=\6R_{00}(Z,\zeta,z',\chi';0)+T\cdot \6P_{00}(Z,\zeta,z',\chi';T),
\end{equation}
with $\6P_{00}(Z,\zeta,z',\chi';T)\in {\mathcal O}(\Delta_{0}^N\times \Delta_{0}^N)[z',\chi'][T]$.
In view of (\ref{trio}), (\ref{col}) and (\ref{blabla}), we obtain
\begin{equation}\Label{ruel}
\6R_{00}\left(Z,\zeta,z',\chi';\rho'_{j_0}((z',\theta(Z^k,f^{s_k}(Z^k); Z,z'))),
(\chi',\1{\theta}(Z^k,f^{s_k}(Z^k);\zeta,\chi'))\right)\equiv 0,
\end{equation}
for $(Z,\zeta)\in \M\cap (\Delta_{Z^k}^N\times
(\Delta_{Z^k}^N)^*)$ and $(z',\chi')\in
\Delta_{g^{s_k}(Z^k)}\times (\Delta_{g^{s_k}(Z^k)})^*$. Setting
$z'=g^{s_k}(Z)$ and $\chi'=\1{g^{s_k}}(\zeta)$ in (\ref{ruel}), we
obtain, in view of (\ref{mmz})
\begin{equation}\Label{85}
\6R_{00}\left(Z,\zeta,g^{s_k}(Z),\1{g^{s_k}}(\zeta);\rho'_{j_0}(f^{s_k}(Z),\1{f^{s_k}}(\zeta))\right)\equiv
0,
\end{equation}
for $(Z,\zeta)$ in some neighborhood of $(Z^k, \1{Z^k})$ in $\M$
and hence, by unique continuation, for all $(Z,\zeta)\in \M \cap
(\Delta^N_{0,s_k}\times (\Delta^N_{0,s_k})^* )$. In view of
(\ref{ji}), (\ref{85}) leads to
\begin{equation}\Label{isis}
\6R_{00}(Z,\zeta,g^{s_k}(Z),\1{g^{s_k}}(\zeta);0)=- \rho'_{j_0}(f^{s_k}(Z),\1{f^{s_k}}(\zeta))
\cdot \6P_{00}(Z,\zeta,g^{s_k}(Z),\1{g^{s_k}}(\zeta);\rho'_{j_0}(f^{s_k}(Z),\1{f^{s_k}}(\zeta))),
\end{equation}
for $(Z,\zeta)$ as above. Since $f(M)\subset M'$, we have the
formal identity $\rho'_{j_0}(f(Z)), \1{f}(\zeta))=0$ for
$(Z,\zeta)\in \M$. Therefore, since $f^{s_k}(Z)$ approximates
$f(Z)$ up to order $s_k\geq k$ at $0$, it follows that
\begin{equation}\Label{pertu}
\rho'_{j_0}(f^{s_k}(Z),\1{f^{s_k}}(\zeta))=O(k),\ (Z,\zeta)\in \M.
\end{equation}
In view of (\ref{isis}), (\ref{pertu}) implies that
\begin{equation}\Label{isis1}
\6R_{00}(Z,\zeta,g^{s_k}(Z),\1{g^{s_k}}(\zeta);0)=O(k)
\end{equation}
for $(Z,\zeta)\in \M$. Since for any $k$, $g^{s_k}(Z)$ approximates
$g(Z)$ up to order $s_k\geq k$ at $0$,
the only possibility for \eqref{isis1} to hold is that
\begin{equation}\Label{vanish}
\6R_{00}(Z,\zeta,g(Z),\1{g}(\zeta);0)\equiv 0,\ (Z,\zeta)\in \M.
\end{equation}
In view of (\ref{nonde}) and (\ref{vanish}), condition (ii) in
Proposition \ref{manger} is satisfied for the components
$g_1(Z),\ldots,g_m(Z)$ of $g(Z)$. By Proposition \ref{manger},
there exists a nontrivial polynomial $\Delta (Z,z')\in \C
\{Z\}[z']$ such $\Delta (Z,g(Z))\equiv 0$. This contradicts the
fact that $g(Z)=(g_1(Z),\ldots,g_m(Z))$ is a transcendence basis
of $\MM_N(f(Z))$ over $\MM_N$. This completes the proof of
Proposition \ref{bourde}.
\end{proof}

\begin{proof}[Proof of Theorem {\rm \ref{straight}}] In view of
Proposition \ref{bourde} (i), we just need to prove
parts (ii) and (iii) of the theorem.
We choose the integer $\kappa_0$ given by Proposition
\ref{bourde} (ii) and define $Z_f$ to be the union of all local
real-analytic irreducible components of $\6Z_f\cap (M\times
\C^{N'})$ at $(0,0)$ that contain the germ of
$\Gamma_{f^{\kappa}}\cap (M\times \C^{N'})$
for some $\kappa \geq \kappa_0$.
The inclusion $Z_f\subset M\times M'$ follows from
the construction of $Z_f$ and Proposition \ref{bourde} (ii).
This shows part (ii) of the theorem. Finally, by setting for any $\kappa\geq
\kappa_0$, $M^{\kappa}:=M\cap \Delta^N_{0,\kappa}$,
part (iii) of the theorem follows
from Propositions \ref{bourde} (ii) and \ref{summa}
and the fact that the subset
$\Sigma^{\kappa}$ is nowhere dense in $\Delta_{0,\kappa}^N$.
The proof of Theorem \ref{straight} is complete.
\end{proof}

\section{Proofs of Theorems \ref{approx} and \ref{ncc0}}\Label{ouf}

\begin{proof}[Proof of Theorem {\rm \ref{approx}}]
Without loss of generality, we may assume that $p$ and $p'$ are
the origin in $\C^N$ and $\C^{N'}$ respectively. In the case where
$M$ is generic, Theorem \ref{approx} is then an immediate
consequence of Theorem \ref{straight} (ii). It remains to consider
the non-generic case. If $M$ is not generic, by using the
intrinsic complexification of $M$, we may assume, after a local
holomorphic change of coordinates near $0$, that $M=\2M\times
\{0\}\subset \C_{z}^{N-r}\times \C_w^r$, for some $1\leq r\leq
N-1$ and some real-analytic generic minimal submanifold $\2M$ (see
e.g.\ \cite{BERbook}). By the generic case treated above, for any
positive integer $k$, there exists a local holomorphic map
$g^k\colon (\C^{N-r},0)\to (\C^{N'},0)$ defined in a neighborhood
of $0$ in $\C^{N-r}$, sending $\2M$ into $M'$ and for which the
Taylor series mapping at $0\in \C^{N-r}$ agrees with $z\mapsto
f(z,0)$ up to order $k$. Let $h^k\colon (\C^N,0)\to (\C^{N'},0)$
be the polynomial mapping obtained by taking the Taylor polynomial
of order $k$ at $0$ of each component of the formal map
$f(z,w)-f(z,0)$. Then by setting for every nonnegative integer
$k$, $f^k(z,w):=g^k(z,0)+h^k(z,w)$, the reader can easily check
that the convergent map $f^k\colon (\C^N,0)\to (\C^{N'},0)$
satisfies all the desired properties. The proof of Theorem
\ref{approx} is complete.
\end{proof}

\begin{proof}[Proof of Theorem {\rm \ref{ncc0}}]
Without loss of generality, we may assume that $p$ and $p'$ are
the origin in $\C^N$ and $\C^{N'}$ respectively. Suppose that $f$
is not convergent. Let $\mu (f)$ and $(f^{\kappa})_{\kappa\geq 0}$
be given by Theorem \ref{straight}. By Proposition \ref{bof}, we
have $m=\mu (f)-N>0$. Therefore Theorem \ref{straight} (iii)
implies that for $\kappa$ large enough, $f^{\kappa}$ maps a dense
subset of a neighborhood of $0$ (which may depend on $\kappa$)
into the subset $\6 E'$. Since $\6 E'$ is closed in $M'$ (see
e.g.\ \cite{D}), $f^{\kappa}$ maps actually a whole neighborhood
of $0$ in $M$ to $\6 E'$. Since for any $\kappa$, $f^{\kappa}$
agrees with $f$ up to order $\kappa$ at $0$, it follows that $f$
sends $M$ into $\6 E'$ as defined in \S\ref{int}. This completes
the proof.
\end{proof}


\begin{thebibliography}{MMZ02a}

\bibitem[A68]{A68} {\bf Artin, M.} --- On the solutions of analytic equations. {\em Invent. math.} {\bf 5} (1968), 277--291.

\bibitem[BER96]{BERacta} {\bf Baouendi,~M.S.; Ebenfelt,~P.; Rothschild, L.P.} --- Algebraicity of holomorphic mappings between real algebraic sets in ${\C}^n$. {\em Acta Math.} {\bf 177} (1996), 225--273.

\bibitem[BER97]{BER97} {\bf Baouendi,~M.S.; Ebenfelt,~P.; Rothschild,~L.P.} --- Parametrization of local biholomorphisms of real analytic hypersurfaces. {\em Asian~J.~Math.} {\bf 1} (1997), 1--16.

\bibitem[BER99a]{BERbook} {\bf Baouendi,~M.S.; Ebenfelt,~P.; Rothschild,~L.P.} --- {\em Real Submanifolds in Complex Space and Their Mappings}. Princeton Math. Series {\bf 47}, Princeton Univ. Press, 1999.

\bibitem[BER99b]{BER99} {\bf Baouendi,~M.S.; Ebenfelt,~P.; Rothschild,~L.P.} --- Rational dependence of smooth and analytic CR mappings on their jets. {\em Math. Ann.} {\bf 315} (1999), 205--249.

\bibitem [BER00a]{BER00} {\bf Baouendi,~M.S.; Ebenfelt,~P.; Rothschild,~L.P.} --- Convergence and finite determination of formal CR mappings. {\em J. Amer. Math. Soc.} {\bf 13} (2000), 697--723.

\bibitem[BER00b]{BERalg}{\bf Baouendi,~M.S.; Ebenfelt,~P.; Rothschild,~L.P.} --- Dynamics of the Segre varieties of a real submanifold in complex space, {\em J. Algebraic Geom.}, to appear.

\bibitem[BG77]{BG} {\bf Bloom,~T.; Graham,~I.} --- On type conditions for generic real submanifolds of $\C^n$. {\em Invent. Math.} {\bf 40} (1977), 217--243.

\bibitem[BMR00]{BMR} {\bf Baouendi, M.S.; Mir, N.; Rothschild, L.P.} --- Reflection ideals and mappings between generic submanifolds in complex space. {\em J. Geom. Anal.}, to appear.

\bibitem[BRZ00]{BRZ} {\bf Baouendi,~M.S.; Rothschild,~L.P.; Zaitsev,~D.} --- Equivalences of real submanifolds in complex space. {\em J. Differential Geom.}, to appear.

\bibitem[CM74]{CM} {\bf Chern,~S.S; Moser,~J.K.} --- Real hypersurfaces in complex manifolds. {\em Acta Math.} {\bf 133} (1974), 219--271.

\bibitem [CPS00]{CPS00} {\bf Coupet,~B.; Pinchuk,~S.; Sukhov,~A.} --- On partial analyticity of CR mappings. {\em Math. Z.} {\bf 235} (2000), 541--557.

\bibitem [D91]{D} {\bf D'Angelo,~J.P.} --- Finite type and the intersection of real and complex subvarieties. {\em Several complex variables and complex geometry, Part 3 (Santa Cruz, CA, 1989)}, 103--117, {\em Proc. Sympos. Pure Math.} {\bf 52}, Part 3, Amer. Math. Soc., Providence, RI, 1991.

\bibitem[K72]{Kohn} {\bf Kohn,~J.J.} --- Boundary behavior of $\bar \partial$ on weakly pseudo-convex manifolds of dimension two. {\em J. Differential Geom.} {\bf 6}, (1972), 523--542.

\bibitem[La01]{la} {\bf Lamel, B.} --- Holomorphic maps of real submanifolds in complex spaces of different dimensions. {\em Pacific J. Math.} {\bf 201} (2001), no. 2, 357--387.

\bibitem[Le86]{Le86} {\bf Lempert,~L.} --- On the boundary behavior of holomorphic mappings. {\em Contributions to several complex variables, 193--215, Aspects Math.}, E9, Viehweg, Braunschweig, 1986.

\bibitem[M00a]{M00} {\bf Mir,~N.} --- Formal biholomorphic maps of real analytic hypersurfaces. {\em Math. Res. Lett.} {\bf 7} (2000), 343--359.

\bibitem[M00b]{M01} {\bf Mir,~N.} --- On the convergence of formal mappings. {\em Comm. Anal. Geom.}, to appear.

\bibitem[MW83]{MW} {\bf Moser, J.K.; Webster, S.M.} --- Normal forms for real surfaces in ${\C}\sp{2}$ near complex tangents and hyperbolic surface transformations. {\em Acta Math.} {\bf 150} (1983), no. 3-4, 255--296.

\bibitem [MMZ02a]{MMZ1} {\bf Meylan,~F.; Mir,~N.; Zaitsev,~D.} --- Analytic regularity of CR-mappings. {\em Math. Res. Lett.} {\bf 9} (2002), 73--93.

\bibitem [MMZ02b]{MMZ01} {\bf Meylan,~F.; Mir,~N.; Zaitsev,~D.} --- Holomorphic extension of smooth CR-map\-pings  between real-analytic and real-algebraic CR-manifolds. Preprint (2002); {\tt http://www.arxiv.org/abs/math.CV/0201267}.

\bibitem[P90a]{Pu1} {\bf Pushnikov,~A.Yu.} --- Holomorphicity of CR-mappings into a space of large dimension. {\em Mat. Zametki} {\bf 48} (1990), no. 3, 147--149.

\bibitem[P90b]{Pu2} {\bf Pushnikov,~A.Yu.} --- On the holomorphy of CR-mappings of real analytic hypersurfaces. {\em Complex analysis and differential equations}, 76--84, Bashkir. Gos. Univ., Ufa, 1990.

\bibitem [S96]{St} {\bf Stanton,~N.} --- Infinitesimal CR automorphisms of real hypersurfaces. {\em Amer. J. Math.} {\bf 118} (1996), 209--233.

\bibitem[T88]{T} {\bf Tumanov, A.E.} --- Extension of CR-functions into a wedge from a manifold of finite type. {\em Mat. Sb. (N.S.)} {\bf 136 (178)} (1988), no. 1, 128--139; translation in {\em Math. USSR-Sb.} {\bf 64} (1989), no. 1, 129--140.

\bibitem[Z97]{Z97} {\bf Zaitsev, D.} --- Germs of local automorphisms of real analytic CR structures  and analytic dependence on the $k$-jets. {\em Math. Res. Lett.} {\bf 4} (1997), 1--20.

\bibitem[Z99]{Z99} {\bf Zaitsev,~D.} --- Algebraicity of local holomorphisms between real-algebraic submanifolds of complex spaces. {\em Acta Math.} {\bf 183} (1999), 273--305.

\bibitem [ZS58]{ZS}  {\bf Zariski,~O.; Samuel,~P.} --- {\em Commutative algebra.} Springer-Verlag, Volume 1, 1958.

\end{thebibliography}
\end{document}